\documentclass[leqno,a4paper,12pt]{article}
\usepackage{amsmath,amsthm,bbm}
\usepackage[sans]{dsfont}

\newtheorem{Thm}{\indent Theorem}[section]
\newtheorem{Prop}[Thm]{\indent Proposition}

\newtheorem{Cor}[Thm]{\indent Corollary}

\theoremstyle{definition}
\newtheorem{Def}[Thm]{\indent Definition}
\newtheorem{Rem}[Thm]{\indent Remark}

\def\qed{{\hskip0pt\unskip\unskip\nobreak\hfil\penalty50
          \hskip1em\hbox{}\nobreak\hfil
          {\bf q.e.d.}%
          \parfillskip=0pt\finalhyphendemerits=0
          \par}\medskip}

\newenvironment{Proof}
               {{\it Proof.}\quad}
               {\qed}

\newenvironment{Proofof}[1]
               {{\it Proof of #1.}\quad}
               {\qed}

\newcommand{\Prime}{\kern3\fontdimen1\font$'$\kern-7\fontdimen1\font}

\long\def\forget#1{}

\long\def\beginSIDEREMARK#1\endSIDEREMARK
    {{\par\bigskip\advance\leftskip by 2cm
                  \advance\rightskip by -2cm\noindent
      {\bf Our own side remark:} #1
      \par\bigskip\noindent}}

\long\def\beginFORGET#1\endFORGET{#1}
\long\def\beginFORGET#1\endFORGET{}

\def\?{\ ???\ \immediate\write16{}%
\immediate\write16{Warning: There was still a question mark . . . }%
\immediate\write16{}}

\usepackage{amsmath}

\usepackage{exscale}
\usepackage{amssymb}

\input cyracc.def
\font\tencyr=wncyr6
\def\cyr{\tencyr\cyracc}
\newcommand{\cyrb}{{\cyr B}}

\newcommand{\B}{{\rm{B}}}

\newcommand{\BA}{{\mathbb{A}}}

\newcommand{\BC}{{\mathbb{C}}}
\newcommand{\BD}{{\mathbb{D}}}

\newcommand{\BF}{{\mathbb{F}}}
\newcommand{\BG}{{\mathbb{G}}}

\newcommand{\BN}{{\mathbb{N}}}

\newcommand{\BQ}{{\mathbb{Q}}}
\newcommand{\BR}{{\mathbb{R}}}

\newcommand{\BZ}{{\mathbb{Z}}}

\newcommand{\Fc}{{\mathfrak{c}}}

\newcommand{\Fe}{{\mathfrak{e}}}

\newcommand{\Fn}{{\mathfrak{n}}}

\newcommand{\FH}{{\mathfrak{H}}}

\newcommand{\FS}{{\mathfrak{S}}}

\newcommand{\FX}{{\mathfrak{X}}}

\newcommand{\CB}{{\cal B}}

\newcommand{\CD}{{\cal D}}

\newcommand{\CF}{{\cal F}}

\newcommand{\CH}{{\cal H}}

\newcommand{\CL}{{\cal L}}

\newcommand{\CV}{{\cal V}}
\newcommand{\CW}{{\cal W}}

\newcommand{\diag}{\mathop{\rm diag}\nolimits}
\newcommand{\Spec}{\mathop{{\bf Spec}}\nolimits}

\newcommand{\GL}{\mathop{\rm GL}\nolimits}
\newcommand{\Gm}{\mathop{\BG_{m,\BQ}}\nolimits}

\newcommand{\Hom}{\mathop{\rm Hom}\nolimits}

\newcommand{\Rep}{\mathop{\rm Rep}\nolimits}

\newcommand{\Sym}{\mathop{\rm Sym}\nolimits}

\newcommand{\loccit}{[loc.$\;$cit.]}

\def\halb{\frac{1}{2}}

\def\id{{\rm id}}

\newbox\mybox
\def\arrover#1{\mathrel{
       \setbox\mybox=\hbox spread 1.4em{\hfil$\scriptstyle#1$\hfil}
       \vbox{\offinterlineskip\copy\mybox
             \hbox to\wd\mybox{\rightarrowfill}}}}
\def\larrover#1{\mathrel{
       \setbox\mybox=\hbox spread 1.4em{\hfil$\scriptstyle#1$\hfil}
       \vbox{\offinterlineskip\copy\mybox
             \hbox to\wd\mybox{\leftarrowfill}}}}

\def\ontoover#1{\mathrel{
       \setbox\mybox=\hbox spread 1.4em{\hfil$\scriptstyle#1$\hfil}
       \vbox{\offinterlineskip\copy\mybox
             \hbox to\wd\mybox{\rightarrowfill\hskip-2.8mm
                               $\rightarrow$}}}}
\def\leftontoover#1{\mathrel{
       \setbox\mybox=\hbox spread 1.4em{\hfil$\scriptstyle#1$\hfil}
       \vbox{\offinterlineskip\copy\mybox
             \hbox to\wd\mybox{$\leftarrow$\hskip-2.8mm
                               \leftarrowfill}}}}
\def\longto{\longrightarrow}
\def\into{\hookrightarrow}
\def\onto{\ontoover{\ }}
\def\longonto{\ontoover{\ }}
\def\isoto{\arrover{\sim}}

\def\longinto{\lhook\joinrel\longrightarrow}

\usepackage[curve,matrix,arrow,cmtip]{xy}
\NoComputerModernTips

\def\myxymessage{\def\messagetext
   {Here an xy-pic diagram was omitted to speed up compilation . . . }
   \immediate\write16{\messagetext}
   \hbox{\bf \messagetext}}
\def\filxymatrix#1{\myxymessage}
\def\filxyarray#1{\myxymessage}

\newdir^{ (}{{}*!/-3pt/\dir^{(}}
\newdir_{ (}{{}*!/-3pt/\dir_{(}}
\newdir^{ )}{{}*!/+3pt/\dir^{)}}
\newdir_{ )}{{}*!/+3pt/\dir_{)}}

\def\rscript#1{\hbox to 0pt{$\scriptstyle#1$\hss}}

\let\oldbullet\bullet
\def\bullet{{\mathchoice{\oldbullet}%
                        {\oldbullet}%
                        {\scriptscriptstyle\oldbullet}%
                        {\oldbullet}}}

\newcommand{\argdot}{{\;\bullet\;}}

\newcommand{\ua}{\mathop{\underline{\alpha}}\nolimits}

\newcommand{\ujast}{\mathop{j_{!*}}\nolimits}
\newcommand{\upjast}{\mathop{j_{!*}'}\nolimits}
\newcommand{\uppjast}{\mathop{j_{!*}''}\nolimits}

\newcommand{\CHMM}{\mathop{CHM(M^K)}\nolimits}
\newcommand{\CHsMM}{\mathop{CHM^s(M^K)}\nolimits}
\newcommand{\CHQQM}{\mathop{CHM(\BQ)_\BQ}\nolimits}

\newcommand{\CHFXM}{\mathop{CHM(X)_\BQ}\nolimits}

\newcommand{\DBcM}{\mathop{DM_{\text{\cyrb},c}}\nolimits}
\newcommand{\DBPAbcM}{\mathop{DM_{\text{\cyrb},c,\Phi}^{Ab}}\nolimits}
\newcommand{\DBcQpPAbM}{\mathop{\DBPAbcM(\partial (M^K)^*)}\nolimits}

\newcommand{\DBcFBsM}{\mathop{\DBcM(\Bs)_F}\nolimits}

\newcommand{\DBcQQM}{\mathop{\DBcM(\BQ)}\nolimits}

\newcommand{\DBcXM}{\mathop{\DBcM(X)}\nolimits}
\newcommand{\DBcFXM}{\mathop{\DBcM(X)}\nolimits}

\newcommand{\DFTBsM}{\mathop{DMT(\Bs)_F}\nolimits}

\newcommand{\Bs}{\mathop{B_\sigma}\nolimits}

\newcommand{\Ss}{\mathop{S_\sigma}\nolimits}
\newcommand{\bSs}{\mathop{\overline{\Ss}}\nolimits}
\newcommand{\St}{\mathop{S_\tau}\nolimits}

\newcommand{\pis}{\mathop{\pi_\sigma}\nolimits}
\newcommand{\pios}{\mathop{\pi'_\sigma}\nolimits}
\newcommand{\pits}{\mathop{\pi''_\sigma}\nolimits}

\newcommand{\one}{\mathds{1}}

\begin{document}

%

\hfuzz=3pt
\overfullrule=10pt                   


\setlength{\abovedisplayskip}{6.0pt plus 3.0pt}
\setlength{\belowdisplayskip}{6.0pt plus 3.0pt}
\setlength{\abovedisplayshortskip}{6.0pt plus 3.0pt}
\setlength{\belowdisplayshortskip}{6.0pt plus 3.0pt}

\setlength{\baselineskip}{13.0pt}
\setlength{\lineskip}{0.0pt}
\setlength{\lineskiplimit}{0.0pt}

%
%

\title{On the intersection motive of certain Shimura varieties:
the case of Siegel threefolds
\forget{
\footnotemark
\footnotetext{To appear in ....}
}
}
\author{\footnotesize by\\ \\
\mbox{\hskip-2cm
\begin{minipage}{6cm} \begin{center} \begin{tabular}{c}
J\"org Wildeshaus \footnote{
Partially supported by the \emph{Agence Nationale de la
Recherche}, project ``R\'egulateurs
et formules explicites''. }\\[0.2cm]
\footnotesize Universit\'e Paris 13\\[-3pt]
\footnotesize Sorbonne Paris Cit\'e \\[-3pt]
\footnotesize LAGA, CNRS (UMR~7539)\\[-3pt]
\footnotesize F-93430 Villetaneuse\\[-3pt]
\footnotesize France\\
{\footnotesize \tt wildesh@math.univ-paris13.fr}
\end{tabular} \end{center} \end{minipage}
\hskip-2cm}
\\[2.5cm]
}
\date{March 22, 2019}
\maketitle
\begin{abstract}
\noindent
In this article, we construct a Hecke-equivariant Chow motive whose realizations equal
intersection cohomology of Siegel 
threefolds with regular algebraic coefficients. 
As a consequence, we are able to define Grothendieck motives
for Siegel modular forms.  \\

\noindent Keywords: Siegel threefolds, weight structures, intersection motive,
motives for Siegel modular forms.

\end{abstract}


\bigskip
\bigskip
\bigskip

\noindent {\footnotesize Math.\ Subj.\ Class.\ (2010) numbers: 
14G35
(11F32, 11F46, 14C25, 14F20, 14F25, 14F30, 14F32, 19E15, 19F27).
}

\eject

\tableofcontents

\bigskip
\vspace*{0.5cm}


%
%

\setcounter{section}{-1}
\section{Introduction}
\label{Intro}



The purpose of this paper is the construction and analysis of the 
\emph{intersection motive} of Kuga--Sato families over a Siegel threefold
relative to its Satake--(Baily--Borel) compactification. 
As in earlier work on Hilbert--Blumenthal varieties \cite{W7},
Picard surfaces \cite{W8}, and more generally, Picard varieties
of arbitrary dimension \cite{C}, the use 
of the formalism of \emph{weight structures} \cite{Bo} proves
to be successful for dealing with a problem, for which explicit
geome\-trical methods seem inefficient. \\

However, Siegel threefolds present a characteristic feature different from the cases
treated so far: the dimension of the boundary of their Satake--(Baily--Borel) compactification
is equal to one. In particular, it is strictly positive. \\

As a consequence, the context of \emph{geometrical motives}, \emph{i.e.}, motives
over a point, is no longer adapted to the problem. Let us explain why. \\

The present construction, as the preceding ones, 
depends on \emph{absence of weights $-1$ and $0$} in the \emph{boundary motive}.
To prove absence of weights, the idea remains, as previously,
to employ \emph{realizations}. But then, realizations need
to detect weights (and therefore, their absence). One may expect
this to be true in ge\-ne\-ral; let us agree to refer to that
principle as \emph{weight conservativity}. To date,
weight conservativity is \emph{proved} for the restriction of the (generic) $\ell$-adic realization
to the category of \emph{motives of Abelian type} of characteristic zero \cite{W11}. \\

However, unless the boundary of the Baily--Borel compactification of a given Shimura variety $M$
is of dimension zero, its boundary motive, as well as the boundary motive of any
Kuga--Sato family $\B$ over $M$, is in general not of Abelian type; this is in any case true 
if $M$ is a Siegel threefold. Concretely, this means that even if the realization
of the boundary motive were proved to avoid weights $-1$ and $0$, we could not formally
conclude that the same is true for the boundary motive itself. \\

This is where \emph{relative motives}, together with the \emph{formalism of six opera\-tions}
enter. Denoting by $j$ the open 
immersion of $M$ into its Baily--Borel compactification $M^*$,
by $i$ its closed complement, and by $\one_M$ the structural motive over $M$,
there is an exact triangle 
\[
i_*i^*j_* \one_M[-1] \longto j_! \one_M \longto j_* \one_M \longto 
i_*i^*j_* \one_M
\]
of motives over $M^*$. The boundary motive of $M$ is isomorphic to the dual of
the direct image of $i_*i^*j_* \one_M$ under the structure morphism of $M^*$.
More generally, the boundary motive of $\B$ is isomorphic to the dual of
the direct image of $i_*i^*j_* \pi_* (\one_\B)$, where $\pi: \B \to M$
denotes the projection of the Kuga--Sato family $\B$ to its base. \\

It is then true that the relative motive $i_*i^*j_* \pi_* (\one_\B)$ over $M^*$
\emph{is} of Abelian type. \\

Whence our strategy of proof.
First, identify the $\ell$-adic realization of $i_*i^*j_* \pi_* (\one_\B)$,
or more generally, of $i_*i^*j_* \CV$, for direct factors $\CV$ of $ \pi_* (\one_\B)$; 
in the cases where weights $0$ and $1$ are
avoided, weight conservativity tells us that $i_*i^*j_* \CV$ itself avoids weights $0$ and $1$.
Second, apply the direct image $a_*$ associated to the structure morphism $a$ of $M^*$.
It is proper, therefore, the functor $a_*$ is \emph{weight exact}. In particular,
if $i_*i^*j_* \CV$ avoids weights $0$ and $1$, then so does $a_* i_*i^*j_* \CV$. 
The corresponding direct
factor of the boundary motive of $\B$ thus avoids weights $-1$ and $0$. \\ 

It may be useful to remark that if $M$ is a Hilbert--Blumenthal or Picard variety,
then there is essentially no difference between $i_*i^*j_* \CV$ and
its direct image under $a$, since the latter is of relative dimension zero on
the boundary of $M^*$. \\

The passage from geometrical motives to relative motives necessitates a certain
number of technical adjustments. For better legibility, we decided to 
separate these from the present text. The result is \cite{W12}; it contains
in parti\-cular the identification of the boundary motive and the dual
of $a_* i_*i^*j_* \pi_* (\one_\B)$ mentioned above. \\

Compared to the cases treated earlier,
another feature of the boundary of Siegel threefolds is new: its canonical stratification 
is not reduced to a single type of strata. Indeed, in the boundary, one finds a closed stratum
of dimension zero, the so-called \emph{Siegel stratum}, and its complement, the so-called
\emph{Klingen stratum}, which is a disjoint union of (open) modular curves. Control
of the weights avoided by the restrictions of the $\ell$-adic realization 
$R_\ell (i^*j_* \pi_* (\one_\B))$ of $i^*j_* \pi_* (\one_\B)$
to the two strata is related to, but does not \emph{a priori} determine
the weights avoided by $R_\ell (i^*j_* \pi_* (\one_\B))$. In fact, the precise relation
is given by a long exact localization sequence. 
Its control it not obvious.
In an earlier attempt, 
we succeeded to identify sufficiently many terms in this sequence,
and (above all) certain morphisms, to prove absence of
weights $0$ and $1$. This approach is technically difficult; moreover, it does not
use the auto-duality property of the coefficients. Indeed, the device dual to the
localization sequence is the co-localization sequence; even when the coefficients
are auto-dual, the two sequences cannot be related. It turns out that both problems
admit the same solution. Namely, the theory of \emph{intermediate extensions}
allows to represent $R_\ell (i^*j_* \pi_* (\one_\B))$ as an extension of two ``halfs'',
one of which dual to the other, and both related to the intermediate extension
$\ujast \pi_* (\one_\B)$. This observation is equally integrated in \cite{W12};
for our purposes, its concrete interest is to divide by two the number
of cohomologi\-cal degrees for which absence of weights has to be tested, and to reduce
the number of morphism in the localization sequence, which need to be identified,
to zero. \\ 

The r\^ole of the intermediate extension is not only technical. It turns out that
the dual of its direct image under $a$ is canonically isomorphic to the \emph{interior
motive}, which according to \cite{W4} can be defined as soon as the boundary
motive avoids weights $-1$ and $0$. This motivates the slight change of terminology
in the title, as compared to the earlier work mentioned above \cite{W7,W8,C}.  \\

Let us now give a more detailed account of the content of the present article. 
Section~\ref{3} contains the statement of our main result, Theorem~\ref{3Main}.
Denote by $GSp_{4,\BQ}$ the group of symplectic similitudes of a fixed four-dimensional
$\BQ$-vector space $V$. As will be recalled, irreducible representations of $GSp_{4,\BQ}$
are indexed by weights $\ua$ depending on three integral parameters: 
$\ua = \alpha(k_1,k_2,r)$. The weight $\ua$ is dominant if and only if 
$k_1 \ge k_2 \ge 0$; it is regular if and only if $k_1 > k_2 > 0$.
Denote by $V_{\ua}$ the irreducible representation
of highest weight $\ua$. According to the main result from \cite{Anc}
(which will  be recalled in Theorem~\ref{3C}),
there is a Chow motive ${}^{\ua} \CV$ over the Siegel threefold $M$,
whose cohomological (Hodge theoretic or $\ell$-adic) realizations
equal the classical \emph{canonical construction} $\mu(V_{\ua})$. Part~(a)
of Theorem~\ref{3Main} then states that $i^* j_* {}^{\ua} \CV$
is of Abelian type. Part~(b) asserts that if $\ua$ is regular,
then $i^* j_* {}^{\ua} \CV$ avoids weights $0$ and $1$.
It has recently become increasingly important to
determine the precise interval containing $[0,1]$ of weights 
avoided by $i^* j_* {}^{\ua} \CV$. 
Theorem~\ref{3Main}~(b) gives a complete answer: putting $k := \min (k_1-k_2, k_2)$,
the motive $i^* j_* {}^{\ua} \CV$ avoids all the weights between
$-k+1$ and $k$, while both weights $-k$ and $k+1$ \emph{do} occur.
Interestingly, this result does not depend on the level 
of the Siegel threefold. We then list the main consequences of 
this result (Corollaries~\ref{3E}, \ref{3F}, \ref{3G}, \ref{3J}, \ref{3I}), 
applying the theory developed in \cite{W12}. \\

Section~\ref{4} is devoted to the proof of Theorem~\ref{3Main}.
As in previous cases, our control of smooth \emph{toroidal compactifications}
of $M$ is sufficiently explicit to verify that, as stated in Theorem~\ref{3Main}~(a),
the motive $i^* j_* {}^{\ua} \CV$ is indeed of Abelian type. Given this result,
and weight conservativity of the restriction of the $\ell$-adic realization
$R_\ell$, part~(b) of Theorem~\ref{3Main} may be checked on the image
of $i^* j_* {}^{\ua} \CV$ under $R_\ell$. Given that ${}^{\ua} \CV$ realizes to give
$\mu(V_{\ua})$, the restriction of $R_\ell (i^* j_* {}^{\ua} \CV)$ to 
the (Siegel and Klingen) strata can be computed following a standard pattern,
employing Pink's and Kostant's Theorems. This computation (Theorem~\ref{4B})
is considerably simplified by results of Lemma's \cite{Lem}. 
It remains to glue the information coming from the strata, in order to 
get control of the weights on the whole boundary.  
The part of Theorem~\ref{3Main}~(b) asserting that
weights $-k$ and $k+1$ occur in $R_\ell (i^* j_* {}^{\ua} \CV)$ 
(Proposition~\ref{4D}) is the single ingredient requiring a proof longer than any other. \\ 

In the final Section~\ref{5}, we give the necessary ingredients to
perform the construction of the Grothendieck motive associated to
a (Siegel) automorphic form with coefficients in an irreducible representation with
regular highest weight (Definition~\ref{5E}). 
This is the analogue for Siegel threefolds of the main result from \cite{Sc}. 
On the level of Galois representations, our definition coincides with
Weissauer's \cite[Thm.~I]{We}. We also recover Urban's result \cite[Thm.~1]{U} 
on characteristic polynomials associated to Frobenii (Corollary~\ref{5G}). \\

Part of this work was done during visits to Caltech's 
Department of Mathe\-matics (Pasadena),
and to the Erwin Schr\"odinger Institute (Vienna). 
I am grateful to both institutions.
I also wish to thank G.~Ancona, J.I.~Burgos Gil, M.~Cavicchi, F.~D\'eglise, F.~Ivorra, F.~Lemma, 
J.~Tilouine and A.~Vezzani for useful discussions and comments,
as well as the referee for her or his observations and suggestions concerning an earlier version
of this article. \\

{\bf Conventions}: We use the triangulated, $\BQ$-linear categories
$\DBcXM$ of \emph{constructible Beilinson motives} over $X$ 
\cite[Def.~15.1.1]{CD},  
indexed by schemes $X$ over $\Spec \BQ$, which are separated and 
of finite type. 
As in \cite{CD}, the symbol $\one_X$
is used to denote the unit for the tensor product in $\DBcFXM$.
We shall employ the full formalism of six operations developed in
\loccit . The reader may choose to consult \cite[Sect.~2]{H} or
\cite[Sect.~1]{W10} for concise presentations of this formalism. \\

Beilinson motives can be endowed with a canonical weight structure,
thanks to the main results from \cite{H}
(see \cite[Prop.~6.5.3]{Bo} for the case
$X = \Spec k$, for a field $k$ of characteristic zero). We refer to it as the \emph{motivic
weight structure}. Following \cite[Def.~1.5]{W10}, the category 
$\CHFXM$ of \emph{Chow motives} over $X$
is defined as the heart $\DBcM(X)_{w = 0}$ of the motivic weight structure
on $\DBcFXM$. \\

A scheme will be said to be \emph{nilregular}
if the underlying reduced scheme is regular in the usual sense.


\bigskip
%
%

\section{Statement of the main result}
\label{3}



In order to state our main result (Theorem~\ref{3Main}),
let us introduce the situation we are going to consider.
The $\BQ$-scheme $M^K$
is a \emph{Siegel threefold}, and the Chow motive ${}^{\ua} \CV$ over $M^K$ is associated  
to a \emph{dominant weight} $\ua = (k_1,k_2,r) \in \BZ^3$,
$k_1 \ge k_2 \ge 0$ 
(see below for the precise normalizations). Denote by $j$ the open 
immersion of $M^K$ into its 
\emph{Satake--(Baily--Borel) compactification} $(M^K)^*$,
and by $i: \partial (M^K)^* \into (M^K)^*$ 
the immersion of the complement of $M^K$ in $(M^K)^*$
(with the reduced scheme structure, say). 
Recall the following.

\begin{Def}[{cmp.~\cite[Def.~2.1~(a)]{W12}}] 
Let $\CHMM_{\BQ,\partial w \ne 0,1}$
denote the full sub-category of $\CHMM_\BQ$
of objects $V$ such that $i^*j_* V$ is without weights $0$ and $1$. 
\end{Def}
 
Theorem~\ref{3Main} implies that in our setting,
the motive ${}^{\ua} \CV \in \CHMM_\BQ$ belongs to $\CHMM_{\BQ,\partial w \ne 0,1}$
if and only if $\ua$ is \emph{regular}: 
$k_1 > k_2 > 0$. More precisely, putting $k := \min (k_1-k_2, k_2)$,
the motive $i^* j_* {}^{\ua} \CV$
is without weights $-k+1,-k+2,\ldots,k$. 
The proof of Theorem~\ref{3Main} 
will be given in Section~\ref{4}. It is an application of
\cite[Thm.~4.4]{W12}; in order to verify the hypotheses of the latter,
we heavily rely on results from \cite{Lem}. \\

Fix a four-dimensional $\BQ$-vector space $V$, together with a $\BQ$-valued
non-degenerate symplectic bilinear form $J$.

\begin{Def} \label{3A}
The group scheme $G$ over $\BQ$ is defined as the group of symplectic similitudes
\[
G := GSp(V,J) \subset GL (V) \; .
\]
\end{Def}
 
Thus, $G$ is reductive, and for any $\BQ$-algebra $R$, the group $G(R)$ equals
\[
\{ g \in \GL (V \otimes_\BQ R) \; , \; 
\exists \, \lambda(g) \in R^* \; , \; 
J(g \argdot,g \argdot) = \lambda(g) \cdot J(\argdot,\argdot) \} \; .
\] 
In particular, the similitude norm $\lambda(g)$ defines a canonical morphism
\[
\lambda: G \longto \Gm \; .
\]

The group $G$ is split over $\BQ$, and its center $Z(G)$ equals
$\Gm \subset \GL(V)$
(inclusion of scalar automorphisms). Maximal $\BQ$-split tori, together with an inclusion
into a Borel sub-group of $G$, are in bijection with symplectic $\BQ$-bases of $V$, 
in which $J$ acquires the $4 \times 4$-matrix
\[
\left( \begin{array}{cc}
0 & I_2 \\
-I_2 & 0
\end{array} \right) \; ,
\]  
equally denoted by $J$. Here as in the sequel, we denote by $I_2$ the 
$2 \times 2$-matrix representing the identity.
Fix one such basis $(e_1,e_2,e_3,e_4)$, use it to identify $G$ with the sub-group $GSp_{4,\BQ}$ of 
$GL_{4,\BQ}$ of matrices $g$
satisfying the relation
\[
 ^{t}\! g J g = \lambda(g) \cdot J \; ,
\]
the maximal split torus with the sub-group $T$ of diagonal matrices
\[
\{ \diag(a,b,a^{-1} q,b^{-1} q) \in GL_{4,\BQ} \} \; ,
\]
and the Borel sub-group with the sub-group of matrices stabilizing the flag of 
totally isotropic sub-spaces $(e_1)_\BQ \subset (e_1,e_2)_\BQ$ of $V$.
We consider triplets $(k_1,k_2,r) \in \BZ^3$ satisfying the congruence relation
\[
r \equiv k_1 + k_2 \!\!\! \mod 2 \; .
\] 
To such a triplet, let us associate the 
(representation-theoretic) weight 
\[
\alpha(k_1,k_2,r) : T \longto \Gm \; , \;
\diag (a,b,a^{-1} q,b^{-1} q) \longmapsto
a^{k_1} b^{k_2} q^{-\frac{r+k_1+k_2}{2}} \; .
\]
Note that
restriction of $\alpha(k_1,k_2,r)$ to $T \cap Sp (V,J)$
corresponds to the projection onto $(k_1,k_2)$. In particular,
the weight $\alpha(k_1,k_2,r)$ is dominant if and only if $k_1 \ge k_2 \ge 0$;
it is regular if and only if $k_1 > k_2 > 0$. Note also that the composition
of $\alpha(k_1,k_2,r)$ with the cocharacter
\[
\Gm \longto T \; , \; 
x \longmapsto \diag (x,x,x,x)
\]
equals 
\[
\Gm \longto \Gm \; , \; x \mapsto x^{-r} \; .
\]
The character $\lambda$ on $T$ equals $\alpha(0,0,-2)$, and
$\det = \lambda^2$.

\begin{Def} \label{3B}
The analytic space $\CH$ is defined as the sub-space of $M_2 (\BC)$ of those
complex $2 \times 2$-matrices, which are symmetrical, and
whose imaginary part is (positive or negative) definite:
\[
\CH := \{ \tau \in M_2 (\BC) \; , \; {}^t \! \tau = \tau \; \text{and} \;
Im (\tau) \, \text{definite} \, \} \; .
\]   
\end{Def}

The group of real points $G(\BR)$ acts on $\CH$
by analytical automorphisms \cite[Ex.~2.7]{P}. In fact,
$(G,\CH)$ are \emph{pure Shimura data} \cite[Def.~2.1]{P}. 
Their \emph{reflex field} \cite[Sect.~11.1]{P}
equals $\BQ$. Given that $Z(G) = \Gm$, 
the Shimura data $(G,\CH)$ satisfy condition 
$(+)$ from \cite[Sect.~5]{W3}. \\
  
Let us now fix additional data:
\begin{enumerate}
\item[(A)] an open compact sub-group $K$ of $G(\BA_f)$
which is neat \cite[Sect.~0.6]{P}, 
\item[(B)] a triplet $(k_1,k_2,r) \in \BZ^3$ 
satisfying the above congruence
\[
r \equiv k_1 + k_2 \!\!\! \mod 2 \; ,
\] 
and in addition,
\[
k_1 \ge k_2 \ge 0 \; .
\]
In other words, the character $\ua := \alpha(k_1,k_2,r)$ is dominant. 
\end{enumerate}

These data are used as follows.
The Shimura variety
$M^K := M^K (G,\CH)$ 
is smooth over $\BQ$. This is the Siegel threefold of level $K$. 
According to \cite[Thm.~11.16]{P}, 
it admits an interpretation as modular space of Abelian surfaces
with additional structures.
In particular, there is a universal family $\B$ 
of Abelian surfaces over $M^K$. \\

The following result holds in the general context of (smooth) Shimura varieties
of $PEL$-type.

\begin{Thm}[{\cite[Thm.~8.6]{Anc}}] \label{3C}
There is a $\BQ$-linear tensor functor
\[
\widetilde{\mu} : \Rep (G) \longto CHM^s \bigl( M^K \bigr)_\BQ 
\]
from the Tannakian category $\Rep (G)$ of algebraic
representations of $G$ in finite dimensional $\BQ$-vector spaces
to the $\BQ$-linear category $CHM^s (M^K)_\BQ$
of \emph{smooth Chow motives} over $M^K$ \cite[Def.~5.16]{L2}.  
It has the following properties.
\begin{enumerate}
\item[(a)] The composition of $\widetilde{\mu}$
with the cohomological
Hodge theoretic reali\-zation is isomorphic to the \emph{canonical construction} 
functor $\mu_{\bf H}$ 
(e.g.\ \cite[Thm.~2.2]{W1}) to the category
of admissible graded-polarizable variations of Hodge structure on $M^K_\BC$.
\item[(b)] The composition of $\widetilde{\mu}$
with the cohomological
$\ell$-adic reali\-zation is isomorphic to the \emph{canonical construction} functor $\mu_\ell$ 
(e.g.\ \cite[Chap.~4]{W1}) to the category
of lisse $\ell$-adic sheaves on $M^K$.
\item[(c)] The functor $\widetilde{\mu}$ commutes with Tate twists.
\item[(d)] The functor $\widetilde{\mu}$
maps the representation $V$
to the dual of the Chow motive $\pi_*^1 \one_{\B}$ over $M^K$.
\end{enumerate} 
\end{Thm}

Here, we denote by $\pi_*^m \one_{\B}$ the 
$m$-th \emph{Chow-K\"unneth component} of the Chow motive 
$\pi_* \one_{\B}$ over $M^K$  \cite[Thm.~3.1]{DM}. 

\medskip

\begin{Proofof}{Theorem~\ref{3C}}
Parts~(a), (c) and (d) are identical to \cite[Thm.~8.6]{Anc}.

As for part~(b), repeat the proof of \loccit , observing that the $\ell$-adic analogue
of \cite[Prop.~8.5]{Anc} holds (the base change to $\BQ_\ell$ of the
sub-group $G_1$ of $G$ coincides with the Lefschetz group).   
\end{Proofof}

\medskip

Given that the representation on $V$ is faithful, it follows that
any object in the image of $\widetilde{\mu}$ 
is isomorphic to a direct sum of
direct factors of Tate twists of the Chow motive 
$\pi_{n_i,*} \one_{\B^{n_i}}$ associated to $\B^{n_i}$, for suitable $n_i \in \BN$,
where $\pi_{n_i}: \B^{n_i} \to M^K$ denotes the 
$n_i$-fold fibre product of $\B$ over $M^K$. 

\begin{Def} \label{3D}
(a)~Denote by $V_{\ua} \in \Rep (G)$ the irreducible representation
of highest weight $\ua$. \\[0.1cm]
(b)~Define ${}^{\ua} \CV \in \CHsMM_\BQ \subset \CHMM_\BQ$ as 
\[
{}^{\ua} \CV := \widetilde{\mu}(V_{\ua}) \; .
\]
\end{Def}

Given that $V_{\ua}$ is of weight $r$, the cohomological realizations of ${}^{\ua} \CV$
equal zero in (classical, \emph{i.e.}, non-perverse) 
degrees $\ne r$, and $\mu_{\bf H}(V_{\ua})$ (in the Hodge theoretic setting)
resp.\ $\mu_\ell(V_{\ua})$ (in the $\ell$-adic setting) in degree $r$. \\ 

Denote by $j: M^K \into (M^K)^*$ the open immersion of $M^K$ into its
Satake--(Baily--Borel) compactification, by $i: \partial (M^K)^* \into (M^K)^*$    
its complement,  and by $\Phi$ the natural stratification of $\partial (M^K)^*$
(the latter will be made explicit in the beginning of Section~\ref{4}). 
Here is our main result.

\begin{Thm} \label{3Main}
(a)~The motive $i^*j_* {}^{\ua} \CV \in \DBcM (\partial (M^K)^*)$ is a
\emph{$\Phi$-con\-structible motive of Abelian type over $\partial (M^K)^*$} (see Definition~\ref{4a}). \\[0.1cm]
(b)~The motive $i^*j_* {}^{\ua} \CV$ is without weights
\[
-k+1, -k+2, \ldots , k ,
\]
where $k := \min (k_1-k_2, k_2)$. Both weights
$-k$ and $k+1$ \emph{do} occur in $i^*j_* {}^{\ua} \CV$. 
In particular, ${}^{\ua} \CV$ belongs to the sub-category 
$\CHMM_{\BQ,\partial w \ne 0,1}$ of $\CHMM_\BQ$
if and only if $\ua$ is regular.
\end{Thm} 

Theorem~\ref{3Main} should be compared to \cite[Thm.~3.5]{W7}, \cite[Thm.~3.8]{W8},
and \cite[Thm.~3.6, Prop.~3.8, Prop.~3.9]{C}
(see also \cite[Rem.~5.8~(b)]{W12}),
which treat the cases of Hilbert--Blumenthal varieties, of Picard surfaces, and of Picard varieties
of arbitrary dimension, respectively. \\

Theorem~\ref{3Main}
will be proved in Section~\ref{4}. For the rest of the present section,  
assume that $k = \min (k_1-k_2, k_2) \ge 1$, \emph{i.e.}, $k_1 > k_2 > 0$.
Given that according to Theorem~\ref{3Main}~(b), the motive ${}^{\ua} \CV$ belongs 
to $\CHMM_{\BQ,\partial w \ne 0,1} \, $, the \emph{intersection motive of 
$M^K$ relative to $(M^K)^*$ with coefficients in ${}^{\ua} \CV$} is at our disposal:
by definition \cite[Def.~3.7]{W12}, it equals 
\[
a_* \ujast {}^{\ua} \CV \in CHM(\BQ)_\BQ \; ,
\]
where $a: (M^K)^* \to \Spec \BQ$ 
is the structure morphism of $(M^K)^*$. 
By abuse of language, let us abbreviate, and refer to $a_* \ujast {}^{\ua} \CV$ as the
\emph{intersection motive with coefficients in ${}^{\ua} \CV$}.
Let us list the main corollaries of Theorem~\ref{3Main}.  

\begin{Cor} \label{3E}
Denote by $a$ and $\tilde{a}$ the structure morphisms of $(M^K)^*$ and  
$M^K$, respectively, and by $m$ the natural transformation $j_! \to j_*$. 
Assume $k_1 > k_2 > 0$, \emph{i.e.}, $k \ge 1$. \\[0.1cm]
(a)~The motive $\tilde{a}_! {}^{\ua} \CV \in \DBcQQM$ is without weights $-k,-k+1, \ldots, -1$,
and the motive $\tilde{a}_* {}^{\ua} \CV \in \DBcQQM$ is without weights $1,2, \ldots, k$. 
More precisely, 
the exact triangles
\[
a_*i_*i^* \ujast {}^{\ua} \CV[-1] \longto \tilde{a}_! {}^{\ua} \CV  \longto a_*\ujast {}^{\ua} \CV \longto
a_*i_*i^* \ujast {}^{\ua} \CV
\]
and
\[
a_*\ujast {}^{\ua} \CV \longto \tilde{a}_*{}^{\ua} \CV \longto a_*i_*i^! \ujast {}^{\ua} \CV[1] \longto 
a_*\ujast {}^{\ua} \CV[1] 
\]
are weight filtrations (of $\tilde{a}_! {}^{\ua} \CV$) avoiding weights $-k,-k+1, \ldots, -1$,
and (of $\tilde{a}_* {}^{\ua} \CV$) avoiding weights $1,2, \ldots, k$, respectively. \\[0.1cm]
(b)~The intersection motive $a_* \ujast {}^{\ua} \CV \in \CHQQM$
behaves functorially with respect to both $\tilde{a}_! {}^{\ua} \CV$ and 
$\tilde{a}_* {}^{\ua} \CV$. In particular, any endomorphism of 
$\tilde{a}_! {}^{\ua} \CV$ or of
$\tilde{a}_* {}^{\ua} \CV$ induces an endomorphism
of $a_* \ujast {}^{\ua} \CV$. \\[0.1cm]
(c)~Let $\tilde{a}_! {}^{\ua} \CV \to N \to \tilde{a}_* {}^{\ua} \CV$ 
be a factorization of the morphism 
$a_* m: \tilde{a}_! {}^{\ua} \CV \to \tilde{a}_* {}^{\ua} \CV$ through a Chow motive
$N \in \CHQQM$. 
Then the intersection motive $a_* \ujast {}^{\ua} \CV$ is canonically identified with
a direct factor of $N$, with a canonical direct complement.
\end{Cor}

\begin{Proof}
Given Theorem~\ref{3Main}, parts (a), (b) and (c) follow from
\cite[Thm.~3.4]{W12}, \cite[Thm.~3.5]{W12} and \cite[Cor.~2.5]{W4}, respectively.
\end{Proof}

The equivariance statement from
Corollary~\ref{3E}~(b) applies in particular
to endomorphisms coming from the \emph{Hecke algebra} 
$\FH(K,G(\BA_f))$
associated to the neat open compact sub-group $K$ of $G(\BA_f)$. 
Recall that by what was said earlier, the relative Chow motive
${}^{\ua} \CV$ is a direct factor of a Tate twist of 
$\pi_{N,*} \one_{\B^N}$,
where $\pi_N: \B^N \to M^K$ denotes the 
$N$-fold fibre product of the universal Abelian scheme $\B$ over $M^K$. 

\begin{Cor} \label{3F}
Assume $k \ge 1$. Then every element of the Hecke algebra $\FH(K,G(\BA_f))$
acts naturally on the intersection motive $a_* \ujast {}^{\ua} \CV$.
\end{Cor}

\begin{Proof}
Let $T \in \FH(K,G(\BA_f))$.
According to Corollary~\ref{3E}~(b), it suffices to show that 
$T$ acts on $\tilde{a}_* {}^{\ua} \CV$. 

To do so, we refer to \cite[pp.~591--592]{W9}. 
\end{Proof}

\begin{Cor} \label{3G}
Assume $k \ge 1$, and
let $\widetilde{\B^N}$ be any smooth compactification of $\B^N$. Then
the intersection motive $a_* \ujast {}^{\ua} \CV$ is
a direct factor of a Tate twist of the Chow motive
$b_* \one_{\widetilde{\B^N}}$ ($b:=$ the structure morphism of the $\BQ$-scheme
$\widetilde{\B^N}$).
\end{Cor}

\begin{Proof}
The motive ${}^{\ua} \CV$ is a direct factor of a Tate twist of 
$\pi_{N,*} \one_{\B^N}$:
\[
{}^{\ua} \CV \longinto \pi_{N,*} \one_{\B^N} (\ell)[2 \ell] \longonto {}^{\ua} \CV \; ,
\]
for a suitable integer $\ell$. The morphism
\[
a_* m : \tilde{a}_! \pi_{N,*} \one_{\B^N} \longto \tilde{a}_* \pi_{N,*} \one_{\B^N}
\]
factors through the Chow motive $b_* \one_{\widetilde{\B^N}}$, hence so does
\[
a_* m : \tilde{a}_! {}^{\ua} \CV \longto \tilde{a}_* {}^{\ua} \CV \; .
\]
Now apply Corollary~\ref{3E}~(c).
\end{Proof}

\begin{Rem}
When $r \ge 0$, then according to \cite[Lemma~4.13]{Anc}, the Chow motive
${}^{\ua} \CV$ is a direct factor of $\pi_{N,*} \one_{\B^N}$ (no Tate twist
needed). 

In this context, let us recall \cite[Cor.~3.10]{W12}: 
the intersection motive $a_* \ujast {}^{\ua} \CV$
is canonically dual to the $e_{\ua}$-part of the \emph{interior motive} of $\B^N$, where
$e_{\ua}$ is the idempotent endomorphism corresponding to the direct factor
${}^{\ua} \CV$ of $\pi_{N,*} \one_{\B^N}$.
\end{Rem}

\begin{Cor} \label{3J}
Assume $k \ge 1$, \emph{i.e.}, that $\ua$ is regular. Then for all $n \in \BZ$, the natural 
maps
\[
H^n \bigl( (M^K)^* (\BC), \ujast \mu_{\bf H}(V_{\ua}) \bigr)
\longto H^n \bigl( M^K (\BC), \mu_{\bf H}(V_{\ua}) \bigr)
\]
(in the Hodge theoretic setting) and
\[
H^n \bigl( (M^K)^* \times_\BQ \bar{\BQ}, \ujast \mu_\ell(V_{\ua}) \bigr)
\longto H^n \bigl( M^K \times_\BQ \bar{\BQ}, \mu_\ell(V_{\ua}) \bigr)
\]
(in the $\ell$-adic setting) are injective. Dually,
\[
H^n_c \bigl( M^K (\BC), \mu_{\bf H}(V_{\ua}) \bigr)
\longto H^n \bigl( (M^K)^* (\BC), \ujast \mu_{\bf H}(V_{\ua}) \bigr)
\]
and
\[
H^n_c \bigl( M^K \times_\BQ \bar{\BQ}, \mu_\ell(V_{\ua}) \bigr)
\longto H^n \bigl( (M^K)^* \times_\BQ \bar{\BQ}, \ujast \mu_\ell(V_{\ua}) \bigr)
\]
are surjective. In other words, the natural maps from intersection cohomology to
cohomology with coefficients in $\mu_{\bf H}(V_{\ua})$, resp.\ $\mu_\ell(V_{\ua})$ 
identify intersection and interior cohomology. 
\end{Cor}

\begin{Proof}
Write ${}^{\ua} \CV$ as a direct factor of  
$\pi_{N,*} \one_{\B^N}  (\ell)[2 \ell]$,
for a suitable integer $\ell$. 
Given Theorem~\ref{3Main}, we may quote \cite[Rem.~3.13~(a), (b)]{W12}
(for $X = (M^K)^*$, $U = M^K$, $C = \B^N$ and $e = e_{\ua}$).
\end{Proof}

As pointed out in \cite[Rem.~3.13~(c)]{W12}, sheaf theoretic considerations alone
suffice to show (without any further reference to geometry)
that Theo\-rem~\ref{3Main} implies Corollary~\ref{3J}. \\

Corollary~\ref{3J} is already known. Indeed, according to \cite[Prop.~1]{MT}, the result generalizes
to Siegel varieties of arbitrary dimension. (However, the proof of \loccit \ is analytic.)

\begin{Rem} \label{3H}
By \cite[Thm.~4.14]{W4}, control of the reduction of \emph{some}
compactification of $\B^N$ implies control of certain properties
of the $\ell$-adic reali\-zation of the intersection motive $a_* \ujast {}^{\ua} \CV$.
According to \cite[Thm.~VI.1.1]{FC}, there exists a smooth compactification
of $\B^N$ having good reduction at each prime number $p$
not dividing the level $n$ of $K$.

\cite[Thm.~4.14]{W4} then yields the following: 
(a)~for all prime numbers $p$ not dividing $n$, the $p$-adic
realization of $a_* \ujast {}^{\ua} \CV$ is crystalline,
(b)~if furthermore $p$ and $\ell$ are different, then
the $\ell$-adic realization of $a_* \ujast {}^{\ua} \CV$
is unramified at $p$. 
\end{Rem}

\begin{Cor} \label{3I}
Assume $k \ge 1$.
Let $p$ be a prime number not dividing the level of $K$. 
Let $\ell$ be different from $p$.
Then the characteristic polynomials of the following coincide: (1)~the action
of Frobenius $\phi$ on the $\phi$-filtered module associated to
the (crystalline) $p$-adic realization of the intersection motive $a_* \ujast {}^{\ua} \CV$,
(2)~the action of a geometrical Frobenius automorphism at
$p$ on the (unramified) $\ell$-adic realization of $a_* \ujast {}^{\ua} \CV$.  
\end{Cor}

\begin{Proof}
Fix a smooth compactification $\widetilde{\B^N}$ of $\B^N$ with good reduction at $p$
\cite[Thm.~VI.1.1]{FC}. Thus, the $\BQ_p$-scheme $\widetilde{\B^N} \times_\BQ \BQ_p$ is the generic fibre
of a smooth and proper scheme $\widetilde{\CB^N}$ over $\BZ_p$. Let us denote by 
$\widetilde{\CB^N_{\BF_p}}$ its special fibre.

The $\phi$-filtered module associated to $p$-adic \'etale cohomology of
$\widetilde{\B^N} \times_\BQ \bar{\BQ}$
is first isomorphic to \emph{Hyodo--Kato} cohomology 
$H^{{}^\bullet}_{HK} ( \widetilde{\B^N} \times_\BQ \bar{\BQ}_p)$ \cite[Sect.~3.2]{Be},
and this isomorphism can be chosen to be motivic in the sense that it commutes with the action of correspondences
in $\widetilde{\B^N} \times_\BQ \widetilde{\B^N}$ \cite[Sect.~4.15]{DN}.
By definition, Hyodo--Kato cohomology is log-crystalline cohomology of a log-smooth model;
in our case, given good reduction, such a model is given by $\widetilde{\CB^N}$ (with divisor equal to zero).
In other words, Hyodo--Kato cohomology equals crystalline cohomology of $\widetilde{\CB^N}$. 
This identification commutes with the action of correspondences in
$\widetilde{\CB^N} \times_{\BZ_p} \widetilde{\CB^N}$. Finally, crystalline cohomology of $\widetilde{\CB^N}$
equals crystalline cohomology of $\widetilde{\CB^N_{\BF_p}}$.

Altogether, the $\phi$-filtered module associated to $p$-adic \'etale cohomology of
$\widetilde{\B^N} \times_\BQ \bar{\BQ}$ is identified with crystalline cohomology of $\widetilde{\CB^N_{\BF_p}}$
in a way compatible with the action of correspondences in
$\widetilde{\CB^N} \times_{\BZ_p} \widetilde{\CB^N}$. Concretely, this means that given a correspondence
$e$ in $\widetilde{\CB^N} \times_{\BZ_p} \widetilde{\CB^N}$, the action of its generic fibre $e_{\BQ_p}$ on 
$p$-adic \'etale cohomology is identified with the action of its special fibre $e_{\BF_p}$ on crystalline cohomology.

For $\ell \ne p$, smooth and proper base change allows us to identify $\ell$-adic cohomology of
$\widetilde{\B^N} \times_\BQ \bar{\BQ}$ and $\ell$-adic cohomology of 
$\widetilde{\CB^N_{\BF_p}} \times_{\BF_p} \bar{\BF}_p$, again compatibly with correspondences.

According to Corollary~\ref{3G}, there is an idempotent endomorphism $e_{\BQ}$ of the Chow motive
associated to $\widetilde{\B^N}$, in other words, an idempotent corres\-pondence in 
$\widetilde{\B^N} \times_\BQ \widetilde{\B^N}$,
whose images in the endomorphism rings of the realizations are projectors onto
the realizations of $a_* \ujast {}^{\ua} \CV$. 
We claim that 
$e_{\BQ_p} := e_\BQ  \times_\BQ \BQ_p$ 
can be extended idempotently to
$\widetilde{\CB^N} \times_{\BZ_p} \widetilde{\CB^N}$. Indeed, accor\-ding to 
\cite[Prop.~5.1.1]{O'S}, the restriction morphism from the
endomorphism ring of the Chow motive associated to
$\widetilde{\CB^N}$ to the one of the Chow motive associated to $\widetilde{\B^N}$
is epimorphic, with nilpotent kernel. 
We now follow a standard line of argument (cmp.\ \cite[proof of Cor.~7.8]{Ki}):
let $\Fe$ be \emph{any} extension
of $e_{\BQ_p}$ to $\widetilde{\CB^N} \times_{\BZ_p} \widetilde{\CB^N}$. The difference
$\Fe - \Fe^2$ is nilpotent, say 
\[
\bigl( \Fe - \Fe^2 \bigr)^N = 0 \; .
\]
But then, 
\[
e_{\BZ_p}:= \bigl( \id_{\widetilde{\CB^N}} - (\id_{\widetilde{\CB^N}} - \Fe)^N \bigr)^N
\]
equally extends $e_{\BQ_p}$ to $\widetilde{\CB^N} \times_{\BZ_p} \widetilde{\CB^N}$.
and $e_{\BZ_p}$ is idempotent. 
Altogether,
there is a smooth and proper scheme $\widetilde{\CB^N_{\BF_p}}$ 
over $\BF_p$, and an idempotent endomorphism 
$e_{\BF_p}$ of the Chow motive associated to $\widetilde{\CB^N_{\BF_p}}$,
whose images in the endomorphism rings of crystalline and $\ell$-adic cohomology,
respectively, are projectors onto the realizations of $a_* \ujast {}^{\ua} \CV$. 
The claim thus follows from
\cite[Thm.~2.~2)]{KM}.
\end{Proof}


\bigskip
%
%

\section{Proof of the main result}
\label{4}



We keep the notation of the preceding section. In order to prove Theorem~\ref{3Main},
the idea is to apply the criterion from \cite[Cor.~4.6]{W12}. \\

In order to check the hypotheses of \loccit,
we first need to fix a finite strati\-fication $\Phi$ of $\partial (M^K)^*$
by locally closed sub-schemes. The canonical
choice would be the restriction $\Phi'$
to $\partial (M^K)^*$ of the natural (finite) stra\-ti\-fication of
$(M^K)^*$ from
\cite[Main Theorem~12.3~(c)]{P}, in other words, all the
strata of $(M^K)^*$ except the open one, \emph{i.e.}, except $M^K$.
According to \cite[Lemma~8.2~(a)]{W9}, $\Phi'$ is \emph{good}, meaning that 
the closure of every stratum is a union of strata.
Furthermore \cite[Lemma~8.2~(b)]{W9}, all strata, 
denoted $i_g (M^{\pi_1(K_1)})$, are smooth over 
$\BQ$ (recall that $K$ is assumed neat, and that $(G,\CH)$ satisfy condition $(+)$),
hence regular. The same is therefore true for the following coarser
stratification $\Phi = \{ 0,1 \}$ of $\partial (M^K)^*$: denote by $i_0: Z_0 \into \partial (M^K)^*$
the disjoint union of all closed strata of $\Phi'$, and by $i_1: Z_1 \into \partial (M^K)^*$
the disjoint union of all strata of $\Phi'$, which are open in $\partial (M^K)^*$.  
Indeed, according to \cite[Sect.~6.3, Ex.~4.25 (with $g=2$)]{P}, 
\[
\partial (M^K)^* = Z_0 \coprod Z_1 \; ;
\]
more precisely, $Z_0$ is of dimension zero, and $Z_1$ of dimension one
(hence so is the whole of $\partial (M^K)^*$). Let us refer
to $Z_0$ as the \emph{Siegel stratum}, and to $Z_1$ as the \emph{Klingen stratum} of
$\partial (M^K)^*$. When it will be necessary to insist on the structure of
stratified scheme of $\partial (M^K)^*$, we shall write $\partial (M^K)^*(\Phi)$
instead of $\partial (M^K)^*$.

\begin{Def}[{\cite[Def.~3.4 and 3.5]{W11}}]  \label{4a}
(a)~Let $S(\FS) = \coprod_{\sigma \in \FS} \Ss$
be a good stratification of a scheme $S(\FS)$.
A morphism $\pi: S(\FS) \to \partial (M^K)^*(\Phi)$ is said to be a
\emph{morphism of good stratifications} if
the pre-image of any of the strata $Z_0$, $Z_1$ of $\partial (M^K)^*$ 
is a union of strata $\Ss \, $. \\[0.1cm]
(b)~A morphism $\pi: S(\FS) \to \partial (M^K)^*(\Phi)$ of good stratifications
is said to be \emph{of Abelian type} if it is proper, and if
the following conditions are satisfied. 
\begin{enumerate}
\item[(1)] All strata $\Ss$, $\sigma \in \FS$,
are nilregular, and for any immersion 
$i_\tau: \St \into \bSs$ of a stratum $\St$ into the closure $\bSs$
of a stratum $\Ss$, the functor $i_\tau^!$
maps $\one_{\bSs}$ to a \emph{Tate motive over $S_\tau$} \cite[Sect.~3.3]{L3}. 
\item[(2)] For all $\sigma \in \FS$
such that $\Ss$ is a stratum
of $\pi^{-1}(Z_m)$, $m \in \{ 0,1 \}$, the morphism
$\pis : \Ss \to Z_m$ can be factorized,
\[
\pis = \pios \circ \pits : \Ss \stackrel{\pits}{\longto} \Bs 
\stackrel{\pios}{\longto} Z_m \; ,
\]
such that the motive 
\[
\pi''_{\sigma,*} \one_{\Ss} \in \DBcFBsM
\]
belongs to the category $\DFTBsM$ of Tate motives over $\Bs \, $,
the morphism $\pios$ is proper and smooth, 
and its pull-back to any geometric
point of $Z_m$ lying over a generic point
is isomorphic to a finite disjoint union of Abelian varieties.
\end{enumerate}  
(c)~An object $V \in \DBcM (\partial (M^K)^*)$ is a 
\emph{$\Phi$-constructible motive of Abelian type over $\partial (M^K)^*$} if the following holds: 
the motive $V$ belongs to the strict, full, dense, $\BQ$-linear 
triangulated sub-category $\DBcQpPAbM$
generated by the images under $\pi_*$ of \emph{$\FS$-constructible
Tate motives over $S(\FS)$} \cite[Def.~3.3]{W11}, where
\[
\pi: S(\FS) \longto \partial (M^K)^*(\Phi)
\]
runs through the morphisms of Abelian type with target
equal to $\partial (M^K)^*(\Phi)$. 
\end{Def}

\begin{Thm} \label{4A}
Let $\ua = \alpha(k_1,k_2,r)$, with
$(k_1,k_2,r) \in \BZ^3$ such that
\[
r \equiv k_1 + k_2 \!\!\! \mod 2 \quad \text{and} \quad k_1 \ge k_2 \ge 0 \; ,
\] 
and consider ${}^{\ua} \CV = \widetilde{\mu}(V_{\ua}) \in \CHMM_\BQ \,$.
The motive $i^*j_* {}^{\ua} \CV$ 
belongs to the full sub-category $\DBcQpPAbM$ of $\DBcM (\partial (M^K)^*) \,$.
In other words, it is a
$\Phi$-constructible motive of Abelian type over $\partial (M^K)^*$.
\end{Thm}

\begin{Proof}
As recalled earlier, the relative Chow motive
${}^{\ua} \CV$ belongs to the strict, full, dense, $\BQ$-linear triangulated sub-category 
\[
\pi_{N,*} DMT \bigl( \B^N \bigr)_\BQ^\natural
\]
of $\DBcM (M^K)$ 
generated by the images under $\pi_{N,*}$ of the category of Tate motives over $\B^N$.
Here as before, $\pi_N: \B^N \to M^K$ denotes the 
$N$-fold fibre product of the universal Abelian scheme $\B$ over $M^K$. 

The latter equals the projection from a \emph{mixed Shimura variety}:
indeed \cite[Ex.~2.7]{P}, 
the representation $V$ of $G$ is of Hodge type $\{ (-1,0) , (0,-1) \}$.
The same is then true for the $r$-th power $V^N$
of $V$. By \cite[Prop.~2.17]{P}, this allows for the construction
of the \emph{unipotent extension} $(P^N,\FX^N)$ of $(G,\CH)$
by $V^N$. 
The pair $(P^N,\FX^N)$ constitute \emph{mixed
Shimura data} \cite[Def.~2.1]{P}. By construction, they come endowed with
a morphism $\pi_N: (P^N,\FX^N) \to (G,\CH)$ of Shimura data,
identifying $(G,\CH)$ with the pure Shimura data
underlying $(P^N,\FX^N)$. In particular, $(P^N,\FX^N)$ also satisfy
condition $(+)$.
Now \cite[Thm.~11.18]{P} there is an open compact neat subgroup $K_N$ of $P^N(\BA_f)$,
whose image under $\pi_N$ equals $K$, 
such that $\B^N$ is identified with the \emph{mixed Shimura variety}
$M^{K_N}:= M^{K_N} (P^N,\FX^N)$, and such
the morphism $M^{K_N} \to M^K$ induced by the morphism $\pi_N$ of Shimura data
is identified with the structure morphism of $\B^N$.

Choose a smooth \emph{toroidal compactification} $M^{K_N} (\FS):= M^{K_N} (P^N,\FX^N,\FS)$
of $M^{K_N}$, associated to a \emph{$K_N$-admissible complete cone 
decomposition} $\FS$ \cite[Sect.~6.4]{P}. Then \cite[proof of Thm.~9.21]{P}, 
modulo a suitable refinement of
$\FS$, the natural stratification of $M^{K_N} (\FS)$, also denoted $\FS$, satisfies the conclusions
of \cite[Lemma~8.1]{W9}, \emph{i.e.}, it is good, and the closures of all strata
are regular. Note that the unique open stratum equals $M^{K_N}$.
According to \cite[Sect.~6.24, Main Theorem~12.4~(b)]{P}, the morphism
$\pi_N : \B^N = M^{K_N} \to M^K$ extends to a proper, surjective
morphism $M^{K_N} (\FS) \to (M^K)^*$, still denoted $\pi_N$. 
From the description given in \cite[Sect.~7.3]{P}, one sees that $\pi_N$ is a morphism
of stratifications.

According to \cite[Cor.~4.10~(b), Rem.~4.7]{W9}, the category 
\[
\pi_{N,*} DMT_\FS \bigl( M^{K_N} (\FS) \bigr)_\BQ^\natural
\]
is obtained by gluing $\pi_{N,*} DMT ( \B^N )_\BQ^\natural$ and
$\pi_{N,*} DMT_\FS ( \pi_N^{-1} ( \partial (M^K)^* ) )_\BQ^\natural \,$.
In particular, 
\[
i^*j_* {}^{\ua} \CV \in 
\pi_{N,*} DMT_\FS \bigl( \pi_N^{-1} \bigl( \partial (M^K)^* \bigr) \bigr)_\BQ^\natural \; .
\]
But $\pi_N$ is of Abelian type \cite[Lemma~8.4]{W9}; therefore,
\[
\pi_{N,*} DMT_\FS \bigl( \pi_N^{-1} \bigl( \partial (M^K)^* \bigr) \bigr)_\BQ^\natural
\subset \DBcQpPAbM \; .
\]
\end{Proof}

Next, we collect information on the restriction of $i^* j_* R_{\ell,M^K}({}^{\ua} \CV)$
to the strata $Z_0$ and $Z_1$. The following
is essentially due to Lemma \cite[Sect.~4]{Lem}.

\begin{Thm} \label{4B}
Let $\ell$ be a prime number. \\[0.1cm]
(a)~For all integers $n \le r+2$, the perverse cohomology sheaf
\[
H^n i_0^* i^* j_* R_{\ell,M^K}({}^{\ua} \CV)
\]
on $Z_0$ is of weights $\le n - (k_1 - k_2)$. The perverse cohomology sheaf
\[
H^{r+2} i_0^* i^* j_* R_{\ell,M^K}({}^{\ua} \CV)
\]
is non-zero, and pure of weight $(r+2) - (k_1 - k_2)$. \\[0.1cm]
(b)~For all integers $n \le r+2$, the perverse cohomology sheaf
\[
H^n i_1^* i^* j_* R_{\ell,M^K}({}^{\ua} \CV)
\]
on $Z_1$ is of weights $\le n - k_2$. The perverse cohomology sheaf
\[
H^{r+2} i_1^* i^* j_* R_{\ell,M^K}({}^{\ua} \CV)
\]
is non-zero, and pure of weight $(r+2) - k_2$. 
\end{Thm}

The proof of Theorem~\ref{4B} will be given after Remark~\ref{4Bd}.
In order to prepare it, recall \cite[Ex.~4.25]{P}
that $Z_0$ and $Z_1$ correspond bijectively to the $G(\BQ)$-conjugacy classes of proper
\emph{rational boundary components} \cite[Sect.~4.11]{P} of $(G,\CH)$. Indeed, 
the group $G(\BQ)$ acts transitively on the set of totally isotropic sub-spaces
of $V$ of a given, strictly positive dimension. \\

We already fixed a basis $(e_1,e_2,e_3,e_4)$ of $V$, in which our symplectic bilinear
form $J$ acquires the $4 \times 4$-matrix 
\[
\left( \begin{array}{cc}
0 & I_2 \\
-I_2 & 0
\end{array} \right) \; ,
\]  
which we equally denoted by $J$. The sub-spaces $V_0'$ and $V_1'$ generated
by $\{ e_1,e_2 \}$ and $\{ e_1 \}$, respectively, are both totally isotropic. \\

Following \cite[Ex.~4.25]{P}, we put $Q_m := Stab_G(V_m')$, $m = 0,1$.
Let $P_m$ denote the normal sub-group of $Q_m$ underlying the 
rational boundary component $(P_m,\FX_m)$ giving rise to $Z_m$ \cite[Sect.~4.7]{P}, and $W_m$ its unipotent radical
(which equals the unipotent radical of $Q_m$).
Then, still according to \cite[Ex.~4.25]{P}, 
$Q_0$ equals 
\[
\biggl\{ \left( \begin{array}{cc}
q \cdot A & A \cdot M \\
0 & {}^t A ^{-1}
\end{array} \right) \; , \; q \in \Gm \, , A \in GL_{2,\BQ} \, , {}^t M = M \biggr\} \; ,
\] 
\[
P_0 = \biggl\{ \left( \begin{array}{cc}
q \cdot I_2 & M \\
0 & I_2
\end{array} \right) \; , \; q \in \Gm \, , {}^t M = M \biggr\} \; ,
\] 
\[
W_0 = \biggl\{ \left( \begin{array}{cc}
I_2 & M \\
0 & I_2
\end{array} \right) \; , \; {}^t M = M \biggr\} \; ,
\] 
while $Q_1$ equals
\[
\Biggl\{ \left( \begin{array}{cccc}
a & aq^{-1} (bu+dw) & v & aq^{-1} (cu+ew) \\
0 & b & w & c \\
0 & 0 & a^{-1}q & 0 \\
0 & d & -u & e
\end{array} \right) \; , \; a \, , be - cd = q \in \Gm \Biggr\} \; ,
\] 
\[
P_1 = \Biggl\{ \left( \begin{array}{cccc}
be - cd & bu+dw & v & cu+ew \\
0 & b & w & c \\
0 & 0 & 1 & 0 \\
0 & d & -u & e
\end{array} \right) \; , \; be - cd \in \Gm \Biggr\} \; ,
\] 
\[
W_1 = \Biggl\{ \left( \begin{array}{cccc}
1 & u & v & w \\
0 & 1 & w & 0 \\
0 & 0 & 1 & 0 \\
0 & 0 & -u & 1
\end{array} \right) \Biggr\} \; .
\] 
Observe that $Q_0 \cap Q_1$ equals the Borel sub-group of $G$ 
stabilizing the flag $V_1' \subset V_0'$, and that both $Q_0$ and $Q_1$ contain the fixed maximal split torus
\[
T = \{ \diag(a,b,a^{-1} q,b^{-1} q) \; , \; a,b,q \in \Gm  \} \; .
\]
In particular, $T$ is canonically identified with a maximal $\BQ$-split torus of the reductive group $Q_m/W_m$,
for $m= 0,1$. Given a (representation-theoretic) weight $\alpha: T \to \Gm$, let us 
denote by $\alpha_m$ the same application, but with $T$ seen as a sub-group of $Q_m/W_m$, $m= 0,1$. \\

Note that
\[
R_{\ell,M^K}({}^{\ua} \CV) = \mu_\ell(V_{\ua})[-r] \; .
\]
Recall that we denote by $\Phi'$ the natural (finite) stra\-ti\-fication of
$(M^K)^*$ from \cite[Main Theorem~12.3~(c)]{P}, which is finer than $\Phi$.
In order to determine the classical cohomology objects 
$R^n i_m^* i^* j_* \mu_\ell(V_{\ua})$, for $m = 0,1$, and $n \in \BZ$, 
one applies the following standard strategy. 
(1)~By Pink's Theorem \cite[Thm.~(5.3.1)]{P2}, the restriction
of $R^n i_m^* i^* j_* \mu_\ell(V_{\ua})$ to any 
individual stratum $Z'$ of $\Phi'$ contributing to $Z_m$ equals 
\[
R^n i_m^* i^* j_* \mu_\ell(V_{\ua})_{ \vert Z' } = 
\bigoplus_{p+q=n} \mu_{\ell,Z'} 
\bigl( H^p \bigl( H_C/K_W , H^q \bigl( Lie(W_m),V_{\ua} \bigr) \bigr)  \bigr)\; .
\]
Here, $H_C/K_W$ is an arithmetic sub-group (depending on $Z'$) of
$C_m/W_m$ \cite[Sect.~(5.2)]{P2}, where
$C_m$ is the identity component of the Zariski closure
of the centralizer in $Q_m(\BQ)$ of the rational boundary component $(P_m,\FX_m)$ 
\cite[Sect.~(3.7)]{P2}, and $\mu_{\ell,Z'}$
is the canonical construction functor to the category of lisse $\ell$-adic sheaves
on $Z'$. (2)~Apply Kostant's Theorem
\cite[Thm.~3.2.3]{Vo}, in order to identify
$H^q ( Lie(W_m),V_{\ua} )$ as a representation of the reductive group
$Q_m/W_m$; this allows in particular
to obtain its weights, and gives potential information concerning cohomology of
$H_C/K_W$ with coeffients in $H^q ( Lie(W_m),V_{\ua} )$. \\

The Hodge theoretic analogue of the above strategy yields the cohomology objects
of $i_m^* i^* j_* \mu_{\bf H}(V_{\ua})_{ \vert Z' } \, $; this was made explicit in 
\cite[Sect.~4]{Lem}. Since steps~(2) of the $\ell$-adic and the Hodge theoretic strategies
are identical, we may use the computations from \loccit \ in our setting. 

\begin{Prop}[{\cite[Sect.~4.3]{Lem}}] \label{4Ba}
Let $\ua = \alpha(k_1,k_2,r)$, with
$(k_1,k_2,r) \in \BZ^3$ such that
\[
r \equiv k_1 + k_2 \!\!\! \mod 2 \quad \text{and} \quad k_1 \ge k_2 \ge 0 \; .
\] 
(a)~For $m = 0,1$, we have 
\[
H^q ( Lie(W_m),V_{\ua} ) = 0
\]
whenever $q < 0$ or $q > 3$. If $0 \le q \le 3$, then the $Q_m/W_m$-representation $H^q 
( Lie(W_1),V_{\ua} )$
is (non-zero and) irreducible. \\[0.1cm]
(b)~The highest (representation-theoretic) weight of $H^q ( Lie(W_0),V_{\ua} )$, $0 \le q \le 3$, is
\[
\alpha_0(k_1,k_2,r) \quad \text{for} \quad q = 0 \; ,
\]
\[
\alpha_0(k_1,-k_2-2,r) \quad \text{for} \quad q = 1 \; ,
\]
\[
\alpha_0(k_2-1,-k_1-3,r) \quad \text{for} \quad q = 2 \; ,
\]
\[
\alpha_0(-k_2-3,-k_1-3,r) \quad \text{for} \quad q = 3 \; .
\]
(c)~The highest (representation-theoretic) weight of $H^q ( Lie(W_1),V_{\ua} )$, $0 \le q \le 3$, is
\[
\alpha_1(k_1,k_2,r) \quad \text{for} \quad q = 0 \; ,
\]
\[
\alpha_1(k_2-1,k_1+1,r) \quad \text{for} \quad q = 1 \; ,
\]
\[
\alpha_1(-k_2-3,k_1+1,r) \quad \text{for} \quad q = 2 \; ,
\]
\[
\alpha_1(-k_1-4,k_2,r) \quad \text{for} \quad q = 3 \; .
\]
\end{Prop}

\begin{Proof}
Note that given our normalization, we have
\[
\alpha(k_1,k_2,r) =  \lambda (k_1,k_2,-r) 
\]
in the notation of \cite[top of p.~87]{Lem}. 

Part~(a) follows from Kostant's Theorem, and from the following fact (see 
\cite[proof of Lemma~4.8]{Lem} and \cite[proof of Lemma~4.10]{Lem}), valid for both $m=0$
and $m=1$: the set of Weyl representatives for $Q_m$ contains no element of length
$<0$ or $>3$, and exactly one element of respective lengths $0$, $1$, $2$ and $3$.

As for part~(c), we refer to \cite[proof of Lemma~4.10]{Lem}.

\cite[proof of Lemma~4.8]{Lem} contains the complete setting
for the application of Kostant's Theorem for $m=0$, but makes it explicit only
for $H^2 ( Lie(W_0),V_{\ua} )$ and $H^3 ( Lie(W_0),V_{\ua} )$.
The reader will have no difficulty to fill in the missing information  
needed for part~(b).
\end{Proof}

Note that both $Q_0/W_0$ and $Q_1/W_1$ are isomorphic to the direct product
$\Gm \times_\BQ GL_{2,\BQ}$. More precisely,
\[
Q_0/W_0 = P_0/W_0 \times_\BQ GL_{2,\BQ} = \Gm \times_\BQ GL_{2,\BQ} \; ,
\]
the identification given by sending the class of a matrix
\[
\left( \begin{array}{cc}
q \cdot A & A \cdot M \\
0 & {}^t A ^{-1}
\end{array} \right) 
\] 
to the pair $(q,A)$, and 
\[
Q_1/W_1 = P_1/W_1 \times_\BQ \Gm = GL_{2,\BQ} \times_\BQ \Gm \; ,
\]
the identification given by sending the class of a matrix
\[
\left( \begin{array}{cccc}
a & aq^{-1} (bu+dw) & v & aq^{-1} (cu+ew) \\
0 & b & w & c \\
0 & 0 & a^{-1}q & 0 \\
0 & d & -u & e
\end{array} \right) 
\] 
to the pair 
\[
\left( \left( \begin{array}{cc}
b & c \\
d & e
\end{array} \right),aq^{-1} \right) \; .
\]
The restriction of the inverse identification to maximal split tori sends 
\[
\left( q, \left( \begin{array}{cc}
x & 0 \\
0 & x^{-1}y
\end{array} \right) \right) \in P_0/W_0 \times_\BQ GL_{2,\BQ}
\]
to
\[
\diag(qx,qx^{-1}y,x^{-1},xy^{-1}) \in T \subset Q_0/W_0
\] 
for $m=0$, and 
\[
\left( \left( \begin{array}{cc}
x & 0 \\
0 & x^{-1}q
\end{array} \right),y \right) \in P_1/W_1 \times_\BQ \Gm
\]
to 
\[
\diag(yq,x,y^{-1},x^{-1}q) \in T \subset Q_1/W_1
\] 
for $m=1$. \\

In the following, the reader will be particularly careful not to confuse two notions
of \emph{weight} associated to representations of reductive groups: 
the highest weights in the sense of representation theory (\emph{e.g.}, those occurring
in Kostant's Theorem), when the representation
is irreducible, and the weights as determined by the action of the weight
cocharacter \cite[Sect.~1.3]{P}, when the group underlies Shimura data.

\begin{Cor} \label{4Bc}
(a)~The $Q_0/W_0$-representations 
$H^q ( Lie(W_0),V_{\ua} )$, $0 \le q \le 2$, 
are (irreducible and) regular, except when $q=0$ and $k_1=k_2$,
in which case $H^0 ( Lie(W_0),V_{\ua} )$ factors through the quotient 
$\Gm \times_\BQ  \Gm$ of the group 
\[
Q_0/W_0 = \Gm \times_\BQ GL_{2,\BQ}
\]
\emph{via} the determinant on the factor $GL_{2,\BQ}$. 
The restriction to $SL_{2,\BQ} \subset GL_{2,\BQ}$ of 
$H^1 ( Lie(W_0),V_{\ua} )$ is of highest (representation-theoretic) weight $k_1 + k_2 + 2$.
The restriction to $P_0/W_0$
of $H^0 ( Lie(W_0),V_{\ua} )$ is of weight $(r+1) - (k_1+k_2)-1$,
and the restriction of $H^1 ( Lie(W_0),V_{\ua} )$ is of weight $(r+2) - (k_1-k_2)$.  
\\[0.1cm]
(b)~The restriction to $P_1/W_1$
of $H^0 ( Lie(W_1),V_{\ua} )$ is of weight $(r+1) - k_1-1$, and
the restriction of $H^1 ( Lie(W_1),V_{\ua} )$ is of weight $(r+2) - k_2-1$.
\end{Cor}

\begin{Proof}
(a): Given the above identifications, the weight $\alpha_0(n_1,n_2,r)$ on $T$
maps 
\[
\left( q, \left( \begin{array}{cc}
x & 0 \\
0 & x^{-1}y
\end{array} \right) \right) \in P_0/W_0 \times_\BQ GL_{2,\BQ}
\]
to 
\[
\alpha_0(n_1,n_2,r) (\diag(qx,qx^{-1}y,x^{-1},xy^{-1}))
= x^{n_1 - n_2} y^{n_2} q^{-\frac{r-n_1-n_2}{2}} \; .
\]
In particular, the restriction of $\alpha_0(n_1,n_2,r)$ to $T \cap SL_{2,\BQ}$
corresponds to the integer $n_1 - n_2$. The first and the second claim thus follow
from Proposition~\ref{4Ba}~(b).

The weight cocharacter $\Gm \to P_0/W_0 = \Gm$ maps $z$ to $z^2$ \cite[Ex.~4.25, Ex.~2.8]{P}.
Its composition with the inclusion into $T$, and with $\alpha_0(n_1,n_2,r)$ 
yields 
\[
\Gm \longto \Gm \; , \; z \longmapsto z^{-r+n_1+n_2} \; .
\]
The third claim thus follows from Proposition~\ref{4Ba}~(b), and from
the normali\-zation of weights of representations \cite[Sect.~1.3]{P}.

(b): The weight cocharacter $\Gm \to P_1/W_1 = GL_{2,\BQ}$ maps $z$ to 
\[
\left( \begin{array}{cc}
z & 0 \\
0 & z
\end{array} \right) 
\]
\cite[Ex.~4.25, Ex.~2.8]{P}. Given the above identifications, 
its composition with the inclusion into $T$ maps $z$ to
$\diag(z^2,z,1,z)$. Further composition 
with $\alpha_1(n_1,n_2,r)$ then yields 
\[
\Gm \longto \Gm \; , \; z \longmapsto z^{-r+n_1} \; .
\]
The claim thus follows from Proposition~\ref{4Ba}~(c).
\end{Proof}

In order to complete the ingredients needed for 
the computation of the $R^n i_m^* i^* j_* \mu_\ell(V_{\ua})$
according to the strategy (1), (2) sketched ealier in this section, observe that
the group $H_C/K_W$ associated to an 
individual stratum $Z'$ of $\Phi'$ contributing to $Z_m$
is a neat arithmetic sub-group of $GL_2(\BQ)$ for $m=0$ \cite[proof of Lemma~4.8]{Lem},
hence of $SL_2(\BQ)$.
In particular, it is of cohomological dimension one.
For $m=1$, the 
group $H_C/K_W$, being a neat arithmetic sub-group of $\BG_m(\BQ)$, is trivial
\cite[proof of Lemma~4.10]{Lem}.

\begin{Rem} \label{4Bd}
When $m=0$, let $V_2$ denote the standard representation of $SL_{2,\BQ}$, and $u \in \BN$.
Then $\Sym^u V_2 \in \Rep (SL_{2,\BQ})$; in fact, $\Sym^u V_2$ is 
the irreductible representation of highest (representation-theoretic) weight $u$. 
Denote by $g$ the genus
of the quotient of the upper half space by $H_C/K_W$, and by $c \ge 1$ the number of its
cusps. (Thus, $c \ge 3$ if $g = 0$ since $H_C/K_W$ is neat.) Then 
$H^1 ( H_C/K_W , \Sym^u V_2 )$ is of dimension
$(u+1)(2g-2+c)$ if $u \ge 1$, and of dimension $2g-1+c$ if $u=0$. In particular,
\[
H^1 \bigl( H_C/K_W , \Sym^u V_2 \bigr) \ne 0 \; , \forall \, u \in \BN \; .
\]
\end{Rem}

\medskip

\begin{Proofof}{Theorem~\ref{4B}}
(a):~According to Corollary~\ref{4Bc}~(a) and Proposition~\ref{4Ba}~(a),
\begin{enumerate}
\item[(o)] $0 \ne H^0 ( Lie(W_0),V_{\ua} )$ is of weight $(r+1) - (k_1 + k_2) - 1$,
\item[(i)] $0 \ne H^1 ( Lie(W_0),V_{\ua} )$ is of weight $(r+2) - (k_1 - k_2)$,
\end{enumerate} 
and $H^q ( Lie(W_0),V_{\ua} ) = 0$ whenever $q < 0$.
The group $H_C/K_W$ associated to a stratum $Z'$ of $Z_0$
is a neat arithmetic sub-group of $SL_2(\BQ)$.
It is therefore of cohomological dimension one,
and admits no non-zero invariants on regular irreducible representations of $Q_0/W_0 = \Gm \times_\BQ GL_{2,\BQ}$.

According to Proposition~\ref{4Ba}~(a) and
Corollary~\ref{4Bc}~(a), $H^q ( Lie(W_0),V_{\ua} )$, $0 \le q \le 2$,
are irreducible as representations of $Q_0/W_0$, and regular unless $q=0$ and $k_1 = k_2$,
in which case $SL_{2,\BQ}$, hence $H_C/K_W$ acts trivially.
Pink's Theorem and \cite[Prop.~(5.5.4)]{P2} then tell us that 
\begin{enumerate}
\item[(o)] $R^0 i_0^* i^* j_* \mu_\ell(V_{\ua})$ is non-zero if and only if
$k_1 = k_2$, in which case it is of weight $r - (k_1 + k_2)$,
\item[(i)] $0 \ne R^1 i_0^* i^* j_* \mu_\ell(V_{\ua})$ is of weight $(r+1) - (k_1 + k_2) - 1$,
\item[(ii)] $0 \ne R^2 i_0^* i^* j_* \mu_\ell(V_{\ua})$ is of weight $(r+2) - (k_1 - k_2)$,
\end{enumerate} 
and that $R^n i_0^* i^* j_* \mu_\ell(V_{\ua}) = 0$ whenever $n < 0$
(for the non-vanishing statements in (i), (ii), see Remark~\ref{4Bd}).
  
The scheme $Z_0$ is of dimension zero; therefore,
\[
H^n i_0^* i^* j_* R_{\ell,M^K}({}^{\ua} \CV) 
= H^{n-r} i_0^* i^* j_*  \mu_\ell(V_{\ua})
= R^{n-r} i_0^* i^* j_*  \mu_\ell(V_{\ua}) \; .
\]
From (o), (i), (ii) and the vanishing of $R^n i_0^* i^* j_* \mu_\ell(V_{\ua}) = 0$ 
for $n < 0$, we conclude that
\begin{enumerate}
\item[(r)] $H^r i_0^* i^* j_* R_{\ell,M^K}({}^{\ua} \CV)$ is zero if $k_1 > k_2$, and
non-zero of weight $r - (k_1 + k_2)$ if $k_1 = k_2$,
\item[(r+1)] $0 \ne H^{r+1} i_0^* i^* j_* R_{\ell,M^K}({}^{\ua} \CV)$ is of weight 
$(r+1) - (k_1 + k_2) - 1$, 
\item[(r+2)] $0 \ne H^{r+2} i_0^* i^* j_* R_{\ell,M^K}({}^{\ua} \CV)$ is of weight 
$(r+2) - (k_1 - k_2)$,
\end{enumerate} 
and that $H^n i_0^* i^* j_* R_{\ell,M^K}({}^{\ua} \CV) = 0$ whenever $n \le r-1$.

(b):~According to Corollary~\ref{4Bc}~(b) and Proposition~\ref{4Ba}~(a),
\begin{enumerate}
\item[(o)] $0 \ne H^0 ( Lie(W_1),V_{\ua} )$ is of weight $(r+1) - k_1-1$,
\item[(i)] $0 \ne H^1 ( Lie(W_1),V_{\ua} )$ is of weight $(r+2) - k_2-1$,
\end{enumerate} 
and $H^q ( Lie(W_1),V_{\ua} ) = 0$ whenever $q < 0$. 
The group $H_C/K_W$ associated to a stratum $Z'$ of $Z_1$ is trivial.
Pink's Theorem and \cite[Lemma~(5.6.6)]{P2} then tell us that 
\begin{enumerate}
\item[(o)] $0 \ne R^0 i_1^* i^* j_* \mu_\ell(V_{\ua})$ is of weight $(r+1) - k_1-1$,
\item[(i)] $0 \ne R^1 i_1^* i^* j_* \mu_\ell(V_{\ua})$ is of weight $(r+2) - k_2-1$,
\end{enumerate} 
and that $R^n i_1^* i^* j_* \mu_\ell(V_{\ua}) = 0$ whenever $n < 0$.
Furthermore, Pink's Theorem tells us that all classical cohomology objects 
$R^n i_1^* i^* j_* \mu_\ell(V_{\ua})$, $n \in \BZ$ are lisse.
The formula 
\[
H^n i_1^* i^* j_* R_{\ell,M^K}({}^{\ua} \CV) 
= H^{n-r} i_1^* i^* j_*  \mu_\ell(V_{\ua})
= \bigr( R^{n-r-1} i_1^* i^* j_*  \mu_\ell(V_{\ua}) \bigr) [1]
\]
is valid: the first equation comes from 
\[
R_{\ell,M^K}({}^{\ua} \CV) = \mu_\ell(V_{\ua})[-r] \; .
\]
As for the second, note that any lisse $\ell$-adic sheaf $\CF$ on a one-dimensional
regular scheme is a perverse sheaf $\CF'$ up to a shift by $-1$:
\[
\CF = \CF'[-1] \quad \text{and} \quad \CF' = \CF[1] \; .
\]
From (o), (i) and the vanishing of $R^n i_1^* i^* j_* \mu_\ell(V_{\ua}) = 0$ 
for $n < 0$, we conclude that
\begin{enumerate}
\item[(r+1)] $0 \ne H^{r+1} i_1^* i^* j_* R_{\ell,M^K}({}^{\ua} \CV)$ is of weight $(r+1) - k_1$,
\item[(r+2)] $0 \ne H^{r+2} i_1^* i^* j_* R_{\ell,M^K}({}^{\ua} \CV)$ is of weight $(r+2) - k_2$,
\end{enumerate} 
and that $H^n i_1^* i^* j_* R_{\ell,M^K}({}^{\ua} \CV) = 0$ whenever $n \le r$.
\end{Proofof}

For the final step of the proof of Theorem~\ref{3Main}, 
the following commutative diagram of immersions
will be useful.
\[
\vcenter{\xymatrix@R-10pt{
        &
        &     
        &
        & 
        &
        Z_1 \ar@{_{ (}->}[llldd]_-{i'} \ar@{_{ (}->}[llddd]^-{i_1}_-\circ \\
        &
        &
        &
        &
        &
        \\ 
        &
        &
        (M^K)^* - Z_0 \ar@{_{ (}->}[d]^-{j''}_-\circ &
        &
        &
        \\
        M^K \ar@{^{ (}->}[rru]^-{j'}_-\circ \ar@{^{ (}->}[rr]_-{j}^-\circ &
        &
        (M^K)^* &
        \partial (M^K)^* \ar@{_{ (}->}[l]^-{i} &
        &
        \\   
        &
        &
        &
        &
        &
        \\   
        &
        &
        &
        &
        & 
        \\
        &
        &
        Z_0 \ar@{^{ (}->}[uuu]_-{i''} \ar@{^{ (}->}[ruuu]_-{i_0} &
        &
        & 
\\}}
\] 
Immersions situated on the same line are complementary to each other 
(exam\-ple: $j''$ and $i''$), the four immersions marked by ``$\circ$''
are open (exam\-ple: $i_1$), and the other four are closed (example: $i'$).  

\begin{Rem} \label{4Cpre}
Denote by $\tau^{t \le \bullet}_{Z_m}$ 
and $\tau^{t \ge \bullet}_{Z_m}$ 
the truncation functors with respect to the perverse $t$-structure on $Z_m$, $m=0,1$. \\[0.1cm]
(a)~The immersions $j'$ and $i'$
being complementary,
\[
(i')^* \upjast \CF' = \tau^{t \le -1}_{Z_1} (i')^* j_*' \CF'
\]
for any perverse sheaf $\CF'$ on $M^K$ \cite[Prop.~1.4.23]{BBD}. \\[0.1cm]
(b)~The intermediate extension is transitive, \emph{i.e.}, 
\[
\ujast = \uppjast \upjast
\]
\cite[Cor.~1.4.24]{BBD}. 
Application of the functor $(i'')^* j_*''$ to the exact triangle
\[
i_*' \tau^{t \ge 0}_{Z_1} (i')^* j_*' [-1] \longto
\upjast \longto
j_*' \longto
i_*' \tau^{t \ge 0}_{Z_1} (i')^* j_*' 
\]
of functors on perverse sheaves on $M^K$ (see (a)) yields the exact triangle
\[
i_0^* i_{1,*} \tau^{t \ge 0}_{Z_1} (i')^* j_*' [-1] \longto
(i'')^* j_*'' \upjast \longto
i_0^* i^* j_* \longto
i_0^* i_{1,*} \tau^{t \ge 0}_{Z_1} (i')^* j_*' 
\]
The immersions $j''$ and $i''$
being complementary, we have as in (a)
\[
(i'')^* \uppjast \CF'' = \tau^{t \le -1}_{Z_0} (i'')^* j_*'' \CF''
\]
for any perverse sheaf $\CF''$ on $(M^K)^* - Z_0$.
It follows that for any perverse sheaf $\CF'$ on $M^K$, there are exact sequences
of perverse cohomology objects
\[
H^{n-1} \bigl( i_0^* i_{1,*} \tau^{t \ge 0}_{Z_1} i_1^* i^* j_* \CF' \bigr) \longto
H^n \bigl( i_0^* i^* \ujast \CF' \bigr) 
\]
\[
\longto
H^n \bigl( i_0^* i^* j_* \CF' \bigr) 
\longto
H^n \bigl( i_0^* i_{1,*} \tau^{t \ge 0}_{Z_1} i_1^* i^* j_* \CF' \bigr) 
\]
for $n \le -1$, while $H^n ( i_0^* i^* \ujast \CF' ) = 0$ for all $n \ge 0$. \\[0.1cm]
(c)~Recall that
$R_{\ell,M^K}({}^{\ua} \CV) = \mu_\ell(V_{\ua})[-r]$; the variety $M^K$
being of dimension three, the complex 
$R_{\ell,M^K}({}^{\ua} \CV)$ is therefore concentrated in perverse degree $r+3$.
According to our conventions, 
$i_1^* i^* \ujast R_{\ell,M^K}({}^{\ua} \CV) = (i')^* \upjast R_{\ell,M^K}({}^{\ua} \CV)$
thus equals
\[
\bigl( (i')^* \upjast ( R_{\ell,M^K}({}^{\ua} \CV) [r+3] ) \bigr) [-(r+3)] \; .
\]
According to (a), we thus have
\[
i_1^* i^* \ujast R_{\ell,M^K}({}^{\ua} \CV) =
\tau^{t \le r+2}_{Z_1} (i')^* j_*' R_{\ell,M^K}({}^{\ua} \CV) =
\tau^{t \le r+2}_{Z_1} i_1^* i^* j_* R_{\ell,M^K}({}^{\ua} \CV) \; .
\]
Similarly, following (b), 
\[
H^n i_0^* i^* \ujast R_{\ell,M^K}({}^{\ua} \CV) = 0 
\]
for all $n \ge r+3$, and 
there are exact sequences
\[
H^{n-1} i_0^* i_{1,*} \tau^{t \ge r+3}_{Z_1} i_1^* i^* j_* R_{\ell,M^K}({}^{\ua} \CV) 
\longto
H^n i_0^* i^* \ujast R_{\ell,M^K}({}^{\ua} \CV)
\]
\[
\longto
H^n i_0^* i^* j_* R_{\ell,M^K}({}^{\ua} \CV)  
\longto 
H^n i_0^* i_{1,*} \tau^{t \ge r+3}_{Z_1} i_1^* i^* j_* R_{\ell,M^K}({}^{\ua} \CV) 
\]
for $n \le r+2$. \\[0.1cm]
(d)~We claim that
\[
H^n i_0^* i_{1,*} \tau^{t \ge r+3}_{Z_1} i_1^* i^* j_* R_{\ell,M^K}({}^{\ua} \CV) = 0
\]
for all $n \le r+1$. Equivalently,
\[
H^n i_0^* i_{1,*} \tau^{t \ge 3}_{Z_1} i_1^* i^* j_* \mu_\ell(V_{\ua}) = 0
\]
for all $n \le 1$. Indeed, by Pink's Theorem, 
the classical cohomology objects of $i_1^* i^* j_* \mu_\ell(V_{\ua})$
are all lisse.
Applying $\tau^{t \ge 3}_{Z_1}$, we thus get a complex concentrated in classical degrees $\ge 2$
(recall that $Z_1$ is of dimension one).
The same is thus true after application of $i_0^* i_{1,*}$ (recall that inverse images
are $t$-exact for the classical $t$-structure). In other words, the complex
\[
i_0^* i_{1,*} \tau^{t \ge 3}_{Z_1} i_1^* i^* j_* \mu_\ell(V_{\ua})
\]
has trivial cohomology (classical or perverse; recall that $Z_0$ is of dimension zero)
in degrees $\le 1$. \\[0.1cm]
(e)~From (c) and (d), we deduce that
\[
H^n i_0^* i^* \ujast R_{\ell,M^K}({}^{\ua} \CV) \isoto 
H^n i_0^* i^* j_* R_{\ell,M^K}({}^{\ua} \CV) 
\]
for $n \le r+1$, and that $H^{r+2} i_0^* i^* \ujast R_{\ell,M^K}({}^{\ua} \CV)$
equals the kernel of
\[
H^{r+2} i_0^* i^* j_* R_{\ell,M^K}({}^{\ua} \CV)  
\longto 
H^{r+2} i_0^* i_{1,*} \tau^{t \ge r+3}_{Z_1} i_1^* i^* j_* R_{\ell,M^K}({}^{\ua} \CV) \; .
\]
\end{Rem}

\begin{Cor} \label{4C}
Let $\ell$ be a prime number. \\[0.1cm]
(a)~For all $n \in \BZ$, 
\[
H^n i_0^* i^* \ujast R_{\ell,M^K}({}^{\ua} \CV)
\]
is of weights $\le n - (k_1 - k_2)$. \\[0.1cm]
(b)~For all $n \in \BZ$, 
\[
H^n i_1^* i^* \ujast R_{\ell,M^K}({}^{\ua} \CV)
\]
is of weights $\le n - k_2$. The perverse cohomology sheaf
\[
H^{r+2} i_1^* i^* \ujast R_{\ell,M^K}({}^{\ua} \CV)
\]
is non-zero, and pure of weight $(r+2) - k_2$. 
\end{Cor}

\begin{Proof}
Part~(a) follows from Remark~\ref{4Cpre}~(c), (e), and from Theorem~\ref{4B}~(a).
 
Part~(b) follows from Remark~\ref{4Cpre}~(c), and from Theorem~\ref{4B}~(b).
\end{Proof}

Corollary~\ref{4C} suffices to prove the part of Theorem~\ref{3Main}~(b) asserting that
regularity of $\ua$ is sufficient for weights $0$ and $1$ to be avoided by $i^*j_* {}^{\ua} \CV$. 
In order to prove that it is necessary, we need the following statement.

\begin{Prop} \label{4D}
Let $\ell$ be a prime number. Then provided that $k_1 \ge 1$, the perverse cohomology sheaf
\[
H^{r+2} i_0^* i^* \ujast R_{\ell,M^K}({}^{\ua} \CV)
\]
is non-zero, and pure of weight $(r+2) - (k_1 - k_2)$.
\end{Prop}

\begin{Proof}
According to Remark~\ref{4Cpre}~(e), 
\[
H^{r+2} i_0^* i^* \ujast R_{\ell,M^K}({}^{\ua} \CV)
\]
equals the kernel of 
\[
ad: H^{r+2} i_0^* i^* j_* R_{\ell,M^K}({}^{\ua} \CV) \longto
H^{r+2}  i_0^* i_{1,*} \tau^{t \ge r+3}_{Z_1} i_1^* i^* j_* R_{\ell,M^K}({}^{\ua} \CV) 
\]
--- in particular, it is pure of weight $(r+2) - (k_1 - k_2)$ (Theorem~\ref{4B}~(a)) ---, 
\emph{i.e.}, it equals the kernel of
\[
H^2 i_0^* i^* j_* \mu_\ell(V_{\ua}) \longto
H^2  i_0^* i_{1,*} \tau^{t \ge3}_{Z_1} i_1^* i^* j_* \mu_\ell(V_{\ua}) \; .
\]
Thanks to Pink's Theorem, the regularity of $H^2 ( Lie(W_0),V_{\ua} )$
as a representation of $Q_0/W_0$ (Corollary~\ref{4Bc}~(a)), and the fact that
the group $H_C/K_W$ is of cohomological dimension one, 
locally on $Z_0$, the (perverse or classical) sheaf
\[
H^2 i_0^* i^* j_* \mu_\ell(V_{\ua}) = R^2 i_0^* i^* j_* \mu_\ell(V_{\ua})
\]
equals 
\[
\mu_{\ell,Z'} 
\bigl( H^1 \bigl( H_C/K_W , H^1 \bigl( Lie(W_0),V_{\ua} \bigr) \bigr)  \bigr)\; ,
\]
for a stratum $Z'$ of $\Phi'$ contributing to $Z_0$.
Furthermore (Corollary~\ref{4Bc}~(a)), the restriction of $H^1 ( Lie(W_0),V_{\ua} )$
to the group $H_C/K_W$ is
isomorphic to the $(k_1 + k_2 +2)$-nd
symmetric power of the standard representation of $SL_{2,\BQ}$.
By Remark~\ref{4Bd}, 
$H^2 i_0^* i^* j_* \mu_\ell(V_{\ua})_{ \vert Z' }$ is therefore of constant rank
$(k_1+k_2+3)(2g-2+c)$, where $g$ denotes the genus of $H_C/K_W$, and $c$ the number of cusps. 

We claim that the restriction to the same $Z'$ of
\[
H^2  i_0^* i_{1,*} \tau^{t \ge3}_{Z_1} i_1^* i^* j_* \mu_\ell(V_{\ua}) 
\]
is of constant rank $c$. Indeed, according to Remark~\ref{4Cpre}~(d),
the classical cohomology objects of $i_1^* i^* j_* \mu_\ell(V_{\ua})$
are all lisse. Therefore, perverse truncation above degree three equals classical
truncation above degree two (recall that $Z_1$ is of dimension one).
The complex
\[
i_0^* i_{1,*} \tau^{t \ge 3}_{Z_1} i_1^* i^* j_* \mu_\ell(V_{\ua})
\]
is concentrated in degrees $\ge 2$, and we get
\[
H^2  i_0^* i_{1,*} \tau^{t \ge3}_{Z_1} i_1^* i^* j_* \mu_\ell(V_{\ua})
= 
R^0  i_0^* i_{1,*} R^2 i_1^* i^* j_* \mu_\ell(V_{\ua}) \; .
\]
Restriction to $Z'$ yields
\[
H^2  i_0^* i_{1,*} \tau^{t \ge3}_{Z_1} i_1^* i^* j_* \mu_\ell(V_{\ua})_{ \vert Z' }
= \bigoplus_{Z''}  
\bigl( R^0  i_0^* i_{1,*} \bigl( R^2 i_1^* i^* j_* \mu_\ell(V_{\ua})_{ \vert Z'' } \bigr) \bigr)_{ \vert Z' } \; ,
\]
where the direct sum is indexed by all strata $Z''$ contributing to $Z_1$, and containing $Z'$ in their closure.
For every such $Z''$, 
\[
R^2 i_1^* i^* j_* \mu_\ell(V_{\ua})_{ \vert Z'' } = 
\mu_{\ell,Z''} 
\bigl( H^2 \bigl( Lie(W_1),V_{\ua} \bigr) \bigr) 
\]
according to Pink's Theorem (since the group $H_C/K_W$ (for $m=1$!) is trivial).

Denote by $j_1: Z_1 \into Z_1^*$ the Baily--Borel compactification,
and by $i_{01}: \partial Z_1^* \into Z_1^*$ its complement. The immersion $i_1: Z_1 \into (M^K)^*$ admits a 
natural extension $\bar{i}_1: Z_1^* \to (M^K)^*$ \cite[Main~Thm.~12.3~(c), Sect.~7.6]{P},
which is finite. The diagram
\[
\vcenter{\xymatrix@R-10pt{
        Z_1^* \ar[d]_{\bar{i}_1} &
        \partial Z_1^* \ar[d]^{\bar{i}_1} \ar@{_{ (}->}[l]_-{i_{01}} \\
        (M^K)^*  &  
        Z_0 \ar@{_{ (}->}[l]_-{i_0} 
\\}}
\]  
is cartesian up to nilpotent elements. Proper base change therefore
yields the formula
\[
R^0  i_0^* i_{1,*} = R^0  \bar{i}_{1,*} i_{0,1}^* j_{1,*} \; .
\]
The functors $ \bar{i}_{1,*}$ and $i_{0,1}^*$ being exact on sheaves, we have
\[
R^0  i_0^* i_{1,*} \bigl( R^2 i_1^* i^* j_* \mu_\ell(V_{\ua})_{ \vert Z'' } \bigr)
= \bar{i}_{1,*} i_{0,1}^* R^0  j_{1,*}  \mu_{\ell,Z''} 
\bigl( H^2 \bigl( Lie(W_1),V_{\ua} \bigr) \bigr)  \; .
\]
According to Proposition~\ref{4Ba}~(a), $H^2 ( Lie(W_1),V_{\ua} )$ is irreducible
as a re\-presentation of $Q_1/W_1$, hence of $GL_{2,\BQ}$.
Yet another application of Pink's Theorem shows that 
\[
i_{0,1}^* R^0  j_{1,*}  \mu_{\ell,Z''} 
\bigl( H^2 \bigl( Lie(W_1),V_{\ua} \bigr) \bigr) 
\]
is of constant rank one on the intersection 
of $\partial Z_1^*$ with the closure of $Z''$ in $(Z_1)^*$.

Our claim on the rank of 
\[
H^2  i_0^* i_{1,*} \tau^{t \ge3}_{Z_1} i_1^* i^* j_* \mu_\ell(V_{\ua})_{ \vert Z' }
= \bar{i}_{1,*} \bigoplus_{Z''}  
\bigl( i_{0,1}^* R^0  j_{1,*}  \mu_{\ell,Z''} 
\bigl( H^2 \bigl( Lie(W_1),V_{\ua} \bigr) \bigr) \bigr)_{ \vert Z' }
\]
is therefore proven as soon as we establish that the number of points in the 
geometrical fibres 
of the morphism $\bar{i}_1: \partial Z_1^* \to Z_0$
above  $Z' \subset Z_0$ equals $c$. 
This verification can be done on the level of $\BC$-valued points,
where the adelic description of the situation is at our disposal.
More precisely,  write $(G_m,\CH_m) := (P_m,\FX_m)/W_m$ \cite[Prop.~2.9]{P}, $m=0,1$,
for the Shimura data contributing to $\partial (M^K)^*$, and 
$Q_{01}$ for the Borel sub-group $Q_0 \cap Q_1$ of $G$. 
According to \cite[Sect.~6.3]{P}, the diagram of $\BC$-valued
points corresponding to the diagram
\[
\vcenter{\xymatrix@R-10pt{
        Z_1^* \ar[d]_{\bar{i}_1} &
        \partial Z_1^* \ar[d]^{\bar{i}_1} \ar@{_{ (}->}[l]_-{i_{01}} \\
        (M^K)^*  &  
        Z_0 \ar@{_{ (}->}[l]_-{i_0} 
\\}}
\] 
equals
\[
\vcenter{\xymatrix@R-10pt{
        Q_1(\BQ) \backslash \bigl( \CH_1^* \times G(\BA_f)/K \bigr) \ar[d]_{\bar{i}_1} &
        Q_{01}(\BQ) \backslash \bigl( \CH_0 \times G(\BA_f)/K \bigr)
                              \ar[d]^{\bar{i}_1} \ar@{_{ (}->}[l]_-{i_{01}} \\
        G(\BQ) \backslash \bigl( \CH^* \times G(\BA_f)/K \bigr)   &  
        Q_0(\BQ) \backslash \bigl( \CH_0 \times G(\BA_f)/K \bigr) \ar@{_{ (}->}[l]_-{i_0} 
\\}}
\] 
where all maps are induced by canonical inclusions of groups and spaces.
Indeed, the full group $Q_m(\BQ)$ (and not only a sub-group of finite index)
stabilizes $\CH_m$, $m=0,1$, and two rational boundary components of $(G_1,\CH_1)$
are conjugate under $G_1(\BQ)$ if and only if they are conjugate under $G(\BQ)$
(explicit computation, or \cite[(iii) of Remark on p.~91]{P}).
The sub-scheme $Z' \subset Z_0$ equals the image of a Shimura variety associated
to $(G_0,\CH_0)$ under a morphism $i_g$ associated to an element $g \in G(\BA_f)$
\cite[Main Theorem~12.3~(c)]{P};
given the adelic description of $i_g$ from \cite[Sect.~6.3]{P},
we see that under the above identification, any $z \in Z'(\BC)$
equals the class $[h_0,p_0g]$ in 
\[
Q_0(\BQ) \backslash \bigl( \CH_0 \times G(\BA_f)/K \bigr) 
\]
of a pair of the form $(h_0,p_0g)$, with $h_0 \in \CH_0$ and $p_0 \in P_0(\BA_f)$.
Put
\[
Q_0^+(\BQ) := \{ q_0 \in Q_0(\BQ) \; , \; \lambda (q_0) > 0 \} \; ;
\]
this group equals the centralizer in $Q_0(\BQ)$ of $h_0$, and indeed,
of the whole of $\CH_0$. Putting
\[
H_C' := Q_0^+(\BQ) \cap p_0gKg^{-1}p_0^{-1} \; ,
\]
we leave it to the reader to verify that the map
\[
Q_{01}(\BQ) \backslash Q_0(\BQ) / H_C' \longto \bar{i}_1^{-1} (z) \; , \;
[q_0] \longmapsto q_0[h_0,p_0g] = [q_0h_0,q_0p_0g]
\]
is well-defined, and bijective. By strong approximation,
\[
W_0(\BQ) \cdot H_C' = Q_0^+(\BQ) \cap W_0(\BA_f) \cdot p_0gKg^{-1}p_0^{-1} \; .
\]
But 
\[
Q_0/W_0 = P_0/W_0 \times_\BQ GL_{2,\BQ} \; ,
\]
meaning that modulo $W_0$, elements in $P_0$ and in $Q_0$ commute with each other.
Thus,
\[
W_0(\BQ) \cdot H_C' = Q_0^+(\BQ) \cap W_0(\BA_f) \cdot gKg^{-1} \; .
\]
The image of $W_0(\BQ) \cdot H_C'$ under the projection $\pi_0 : Q_0 \onto Q_0/W_0$
coincides with the image of 
\[
W_0(\BA_f) \cdot Q_0^+(\BQ) \cap gKg^{-1} 
\]
(both images equals $\pi_0 (Q_0^+(\BQ)) \cap \pi (gKg^{-1})$).
But by definition \cite[(3.7.4)]{P2}, $W_0(\BA_f) \cdot Q_0^+(\BQ) \cap gKg^{-1}$
equals $H_C$. We thus showed that
\[
\pi_0 (H_C') = \pi_0(H_C) \; .
\]
Now the quotient morphism $Q_0 \onto Q_0/P_0, q_0 \mapsto 
\overline{q_0}$
induces an isomorphism
\[
Q_{01}(\BQ) \backslash Q_0(\BQ) / H_C' 
\isoto \overline{Q_{01}(\BQ)} \backslash \overline{Q_0(\BQ)} / \overline{H_C'}
= \overline{Q_{01}(\BQ)} \backslash \overline{Q_0(\BQ)} / \overline{H_C}
\]
But $\overline{Q_0(\BQ)} = GL_2(\BQ)$, and under this identification,
$\overline{Q_{01}(\BQ)}$ equals the sub-group of upper triangular matrices,
while $\overline{H_C} = H_C/K_W$. In other words, 
\[
Q_{01}(\BQ) \backslash Q_0(\BQ) / H_C'  
\]
is identified with the set up cusps of $H_C/K_W$.

The formula
\[
(k_1 + k_2 + 3)(2g-2+c)  \ge 4(2g-2+c) > c
\]
(recall that $c$ is greater or equal to $1$, and that $c \ge 3$ if $g = 0$)
implies that the rank of the source of $ad$ is strictly greater than the rank of its target; the
kernel of $ad$ is therefore non-trivial.
\end{Proof}

\begin{Rem}
(a)~As the reader may verify, 
\[
H^{r+2}  i_0^* i_{1,*} \tau^{t \ge r+3}_{Z_1} i_1^* i^* j_* R_{\ell,M^K}({}^{\ua} \CV) 
\]
is pure of weight $(r+2) - (k_1 - k_2)$, \emph{i.e.}, of the same
weight as 
\[
H^{r+2} i_0^* i^* j_* R_{\ell,M^K}({}^{\ua} \CV) \; .
\] 
Weight considerations alone do therefore not imply non-triviality of the kernel
of the map $ad$ from the proof of Proposition~\ref{4D}. \\[0.1cm]
(b)~A more conceptual proof 
of Proposition~\ref{4D} would consist in showing that locally on $Z_0$, the map $ad$
equals the direct sum over all cusps of $H_C/K_W$ of the residue maps. 
Identify 
$H^1 ( H_C/K_W , H^1 ( Lie(W_0),V_{\ua} ) ) \otimes_{\BQ} \BC$ with the direct sum of
the space of modular forms and (the conjugate of) the space of cusp forms for $H_C/K_W$ of weight 
$k_1+k_2+4 \ge 5$.
The kernel of the residues contains the space of cusp forms. Its dimension is computed in 
\cite[Thm.~2.24 and Thm.~2.25]{Sh}; thanks to \cite[Prop.~1.40]{Sh} (always remember that
$H_C/K_W$ is neat), this dimension can be seen to be strictly positive. \\[0.1cm]
(c)~On the level of geometry of Baily--Borel compactifications, a ``strange duality'' 
seems to be involved
in the proof of Proposition~\ref{4D}: we need to know how many modular curves in the boundary of
$(M^K)^*$ contain a given cusp $Z'$ in their closure. The response yields the number of cusps of a
``modular curve'', which does not explicitely occur in $(M^K)^*$, namely the quotient of the upper half
space by $H_C/K_W$. 
It would be interesting to see how this phenomenon generalizes to higher dimensional Siegel varieties. \\[0.1cm]
(d)~Our computation of the fibres of the morphism $\bar{i}_1: Z_1^* \to (M^K)^*$ over points of $Z_0$
is a quantitative version of a classical non-injectivity result of Satake \cite[Exemple on p.~13-06]{Sat}.
\end{Rem}

\begin{Rem} 
The Hodge theoretic analogues of Theorem~\ref{4B}, Corollary~\ref{4C} and 
Proposition~\ref{4D} hold. The proofs are identical up to the use of
Pink's Theorem, which is replaced by \cite[Thm.~2.9]{BW}.
\end{Rem}

\medskip

\begin{Proofof}{Theorem~\ref{3Main}}
According to Theorem~\ref{4A}, $i^*j_* {}^{\ua} \CV$ is a
$\Phi$-constructible motive of Abelian type over $\partial (M^K)^*$; this proves part~(a)
of our claim.

By \cite[Summ.~1.18~(d)]{P}, there is a perfect pairing
\[
V_{\ua} \otimes_\BQ V_{\ua} \longto \BQ(-r)
\]
in $\Rep(G)$. 

Fix a prime $\ell$. Applying $\mu_\ell$, we get a perfect pairing  
\[
\mu_\ell(V_{\ua}) \otimes_{\BQ_\ell} \mu_\ell(V_{\ua}) \longto \BQ_\ell(-r)
\]
of $\ell$-adic lisse sheaves on $M^K$. In terms of local duality, the pairing induces
an isomorphism
\[
\BD_{\ell,M^K} \bigl( \mu_\ell(V_{\ua}) \bigr) \cong \mu_\ell(V_{\ua}) (r+3)[6] 
\]   
($M^K$ is smooth of dimension three). Given 
$R_{\ell,M^K}({}^{\ua} \CV) = \mu_\ell(V_{\ua})[-r]$, we find that
\[
\BD_{\ell,M^K} \bigl( R_{\ell,M^K}({}^{\ua} \CV) \bigr) \cong R_{\ell,M^K}({}^{\ua} \CV) (s)[2s] \; , 
\]   
where $s = r+3$. 

Corollary~\ref{4C} tells us that for all $n \in \BZ$, and $m = 0,1$, 
\[
H^n i_m^* i^* \ujast R_{\ell,M^K}({}^{\ua} \CV)
\]
is of weights $\le n - k$. According to \cite[Cor.~4.6~(b)]{W12}, the motive $i^*j_* {}^{\ua} \CV$ 
therefore avoids weights $-k+1, -k+2, \ldots , k$.

In order to conclude the proof of part~(b),
it remains to show, again thanks to \cite[Cor.~4.6~(b)]{W12}, that for some $n \in \BZ$, and 
$m = 0$ or $m=1$, weight $n - k$ \emph{does} occur in
\[
H^n i_m^* i^* \ujast R_{\ell,M^K}({}^{\ua} \CV) \; .
\]
We take $n = r+2$, and distinguish two cases. If $k = k_2$, \emph{i.e.}, $k_2 \le k_1 - k_2$,
take $m=1$; the claim then follows from Corollary~\ref{4C}~(b).  
Else, $k_2 > k_1 - k_2$ and $k = k_1 - k_2$. Since $k_1 \ge k_2$, 
we necessarily have $k_1 \ge 1$. Take $m=0$ and apply Proposition~\ref{4D}.
\end{Proofof}

\medskip

\begin{Rem}
(a)~An element of $H^n ( \partial (M^K)^* (\BC), i^*j_* \mu_{\bf H}(V_{\ua}) )$
is called a \emph{ghost class} if it lies in the image of
\[
H^n \bigl( M^K (\BC), \mu_{\bf H}(V_{\ua}) \bigr) \longto 
H^n \bigl( \partial (M^K)^* (\BC), i^*j_* \mu_{\bf H}(V_{\ua}) \bigr)
\]
and in the kernel of both restriction maps
\[
H^n \bigl( \partial (M^K)^* (\BC), i^*j_* \mu_{\bf H}(V_{\ua}) \bigr) \longto
H^n \bigl( Z_m (\BC), i_m^*i^*j_* \mu_{\bf H}(V_{\ua}) \bigr) \; ,
\]
$m = 0, 1$. One of the main results of \cite{M} implies that if $\ua$ is regular, then 
there are no non-zero ghost classes \cite[Thm.~3.1]{M}. This result does not formally imply,
nor is it implied by our Theorem~\ref{3Main}. Nonetheless, it might be worthwhile to note 
that the weight arguments that occur in the proofs are quite similar. The most relevant information from 
Theorem~\ref{3Main}, as far as \cite[Thm.~3.1]{M} is concerned, comes from the weight
filtration
\[
a_*\ujast {}^{\ua} \CV \longto \tilde{a}_*{}^{\ua} \CV \longto a_*i_*i^! \ujast {}^{\ua} \CV[1] \longto 
a_*\ujast {}^{\ua} \CV[1] 
\]
avoiding weights $1,2, \ldots, k$ (Corollary~\ref{3E}~(a)), hence avoiding weight $1$ if $\ua$ is regular,
which we assume in the sequel.
This implies 
that any element of $H^n ( M^K (\BC), \mu_{\bf H}(V_{\ua}) )$ not mapping to zero in
$H^n ( \partial (M^K)^* (\BC), i^*j_* \mu_{\bf H}(V_{\ua}) )$, remains non-zero in
\[
H^n \bigl( \partial (M^K)^* (\BC), i^! \ujast \mu_{\bf H}(V_{\ua})[1] \bigr) =
H^n \bigl( \partial (M^K)^* (\BC), \tau_{\partial (M^K)^*}^{t \ge 3} i^*j_* \mu_{\bf H}(V_{\ua}) \bigr) \; .
\]
In other words, a ghost class vanishing in 
$H^n ( \partial (M^K)^* (\BC), \tau_{\partial (M^K)^*}^{t \ge 3} i^*j_* \mu_{\bf H}(V_{\ua}) )$ 
is zero. The Hodge structure 
$H^n ( \partial (M^K)^* (\BC), \tau_{\partial (M^K)^*}^{t \ge 3} i^*j_* \mu_{\bf H}(V_{\ua}) )$
has weights $\ge (r+n) +2$; the same type of considerations as those leading to Corollary~\ref{4C}
then imply that the direct sum of the 
restriction maps
\[
H^n \bigl( \partial (M^K)^* (\BC), \tau_{\partial (M^K)^*}^{t \ge 3} i^*j_* \mu_{\bf H}(V_{\ua}) \bigr) \longto
H^n \bigl( Z_m (\BC), i_m^* \tau_{\partial (M^K)^*}^{t \ge 3}i^*j_* \mu_{\bf H}(V_{\ua}) \bigr) \; ,
\]
$m = 0, 1$, is injective. \\[0.1cm]
(b)~The above illustrates an observation made by Moya Giusti: {\it for a class in the cohomology
of the boundary whose weight is neither 
the middle weight nor the middle weight plus one we can determine
exactly whether it is or not in the image of the morphism}
\[
H^n \bigl( M^K (\BC), \mu_{\bf H}(V_{\ua}) \bigr) \longto 
H^n \bigl( \partial (M^K)^* (\BC), i^*j_* \mu_{\bf H}(V_{\ua}) \bigr) \; .
\]
In fact, it appears amusing to note that the ``middle weight'' is relevant in another context 
than the one studied in the present paper. 
According to \cite[p.~2317, second paragraph]{M}, the representation $V_{\ua}$ 
\emph{satisfies the middle weight property} if the space of
ghost classes in $H^n ( \partial (M^K)^* (\BC), i^*j_* \mu_{\bf H}(V_{\ua}) )$ is pure of weight
$r+n$. \cite[Thm.~3.1]{M} implies in particular that for all $\ua$ (regular or not),
the representation $V_{\ua}$ does satisfy the middle weight property, while our Theorem~\ref{3Main}
implies that weights $\{ r+n, r+n+1 \}$ do not occur at all in 
$H^n ( \partial (M^K)^* (\BC), i^*j_* \mu_{\bf H}(V_{\ua}) )$, as soon as $\ua$ is regular.
\end{Rem}

\begin{Rem} \label{4X}
Saper's vanishing theorem \cite[Thm.~5]{Sa} says that if $\ua$ is regular,
then the groups $H^n ( M^K (\BC), \mu_{\bf H}(V_{\ua}) )$,
hence (by comparison) $H^n ( M^K \times_\BQ \bar{\BQ}, \mu_\ell(V_{\ua}) )$
vanish for $n < 3 = \dim M^K$. By duality, one obtains that $H^n_c ( M^K (\BC), \mu_{\bf H}(V_{\ua}) )=0$
and $H^n_c ( M^K \times_\BQ \bar{\BQ}, \mu_\ell(V_{\ua}) )=0$ for $n > 3$. 
It follows that interior cohomology 
with coefficients in $\mu_{\bf H}(V_{\ua})$, denoted
\[
H^n_! \bigl( M^K (\BC), \mu_{\bf H}(V_{\ua}) \bigr) \; ,
\]
and interior cohomology with coefficients
in $\mu_\ell(V_{\ua})$, denoted
\[
H^n_! \bigl( M^K \times_\BQ \bar{\BQ}, \mu_\ell(V_{\ua}) \bigr) \; ,
\]
both vanish for $n \ne 3$, provided that $\ua$ is regular.
\end{Rem}


\bigskip
%
%

\section{The motive for an automorphic form}
\label{5}



This final section contains the analogues for Siegel threefolds
of the main results from \cite{Sc}. Since we shall not restrict
ourselves to the case of Hecke eigenforms, our
notation becomes a little more technical than in \loccit . \\

We continue to consider the situation 
of Sections~\ref{3} and \ref{4}. In particular,
we fix a dominant $\ua = \alpha(k_1,k_2,r)$, which we
assume to be regular, \emph{i.e.}, $k_1 > k_2 > 0$. Consider 
the intersection motive $a_* \ujast {}^{\ua} \CV \in \CHQQM$,
where as before $a: (M^K)^* \to \Spec \BQ$ 
denotes the structure morphism of $(M^K)^*$. According to 
\cite[Rem.~3.13~(a)]{W12} and Remark~\ref{4X},
its Hodge theoretic realization equals
\[
H^3_! \bigl( M^K (\BC), \mu_{\bf H}(V_{\ua}) \bigr)[-(r+3)]
\; ,
\]
and its $\ell$-adic realization equals
\[
H^3_! \bigl( M^K \times_\BQ \bar{\BQ}, \mu_\ell(V_{\ua}) \bigr)[-(r+3)]
\; .
\]
By Corollary~\ref{3F}, every element of the Hecke algebra $\FH(K,G(\BA_f))$ acts on
$a_* \ujast {}^{\ua} \CV$. 

\begin{Thm}[{\cite[Thm.~3.1.1]{Ha}}]  \label{5A}
Let $L$ be any field of characteristic zero. Then the 
$\FH(K,G(\BA_f)) \otimes_\BQ L$-module 
$H^3_! ( M^K (\BC), \mu_{\bf H}(V_{\ua}) ) \otimes_\BQ L$
is semi-simple. 
\end{Thm}

Note that \cite[Sect.~8.1.6, p.~232]{Ha} gives a proof of
Theorem~\ref{5A}, while the statement in \cite[Thm.~3.1.1]{Ha}
is ``non-adelic''. Denote by $R(\FH) := R(\FH(K,G(\BA_f)))$
the image of the Hecke algebra in the endomorphism algebra 
of $H^3_! ( M^K (\BC), \mu_{\bf H}(V_{\ua}) )$. 

\begin{Cor} \label{5B}
Let $L$ be any field of characteristic zero. Then the $L$-algebra
$R(\FH) \otimes_\BQ L$ is semi-simple.
\end{Cor}

In particular, the isomorphism classes of simple right 
$R(\FH)\otimes_\BQ L$-modules correspond bijectively to 
isomorphism classes of minimal right ideals. \\

Fix $L$, and let $Y_{\pi_f}$ be such a minimal right ideal of $R(\FH)\otimes_\BQ L$. 
There is a (primitive) idempotent $e_{\pi_f} \in R(\FH)\otimes_\BQ L$
generating $Y_{\pi_f}$. 

\begin{Def}  \label{5C}
(a)~The \emph{Hodge structure $W(\pi_f)$ associated to $Y_{\pi_f}$} is defined as
\[
W(\pi_f) := \Hom_{R(\FH) \otimes_\BQ L} \bigl( Y_{\pi_f} , 
H^3_! \bigl( M^K (\BC), \mu_{\bf H}(V_{\ua}) \bigr) \otimes_\BQ L \bigr) \; .
\]
(b)~Let $\ell$ be a prime number.
The \emph{Galois module $W(\pi_f)_\ell$ associated to $Y_{\pi_f}$} is defined as
\[
W(\pi_f)_\ell := \Hom_{R(\FH) \otimes_\BQ L} \bigl( Y_{\pi_f} , 
H^3_! \bigl( M^K \times_\BQ \bar{\BQ}, \mu_\ell(V_{\ua}) \bigr) \otimes_\BQ L \bigr) \; .
\]
\end{Def}

Definition~\ref{5C}~(b) should be compared to \cite[Thm.~I]{We}. 

\begin{Prop} \label{5D}
There is canonical isomorphism of Hodge structures
\[
W(\pi_f) \isoto 
\bigl( H^3_! \bigl( M^K (\BC), \mu_{\bf H}(V_{\ua}) \bigr) \otimes_\BQ L \bigr) 
\cdot e_{\pi_f} \; ,
\]
and a canonical isomorphism of Galois modules
\[
W(\pi_f)_\ell \isoto 
\bigl( H^3_! \bigl( M^K \times_\BQ \bar{\BQ}, \mu_\ell(V_{\ua}) \bigr) \otimes_\BQ L \bigr) 
\cdot e_{\pi_f} \; .
\]
\end{Prop}

\begin{Proof}
We shall perform the proof for Hodge structures; the one for Galois modules
is formally identical.
Obviously, 
\[
\Hom_{R(\FH) \otimes_\BQ L} \bigl( R(\FH) \otimes_\BQ L , 
H^3_! \bigl( M^K (\BC), \mu_{\bf H}(V_{\ua}) \bigr) \otimes_\BQ L \bigr)
\]
is canonically identified with
\[
H^3_! \bigl( M^K (\BC), \mu_{\bf H}(V_{\ua}) \bigr) \otimes_\BQ L 
\]
by mapping an morphism $g$ to the image of $1 = 1_{R(\FH)}$ under $g$.
Inside 
\[
\Hom_{R(\FH) \otimes_\BQ L} \bigl( R(\FH) \otimes_\BQ L , 
H^3_! \bigl( M^K (\BC), \mu_{\bf H}(V_{\ua}) \bigr) \otimes_\BQ L \bigr) \; ,
\]
the object $W(\pi_f)$ contains precisely those morphisms $g$ 
vanishing on $1 - e_{\pi_f}$, in other words, satisfying
the relation
$g(1) = g(e_{\pi_f}) = g(1) \cdot e_{\pi_f}$.
\end{Proof}

Since we do not know whether the Chow motive $a_* \ujast {}^{\ua} \CV$
is finite dimensional, we cannot apply \cite[Cor.~7.8]{Ki}, and
therefore do not know whether $e_{\pi_f}$ can be lifted \emph{idempotently}
to the Hecke algebra $\FH(K,G(\BA_f))$. This is why we need to descend
to the level of \emph{Grothendieck motives}. Denote by 
$a_* \ujast {}^{\ua} \CV'$ the Grothendieck motive underlying
$a_* \ujast {}^{\ua} \CV$.

\begin{Def} \label{5E}
Assume $\ua = \alpha(k_1,k_2,r)$ to be regular. Let $L$ be a field
of characteristic zero, and $Y_{\pi_f}$ a minimal right ideal
of $R(\FH)\otimes_\BQ L$.
The \emph{motive associated to $Y_{\pi_f}$} is defined as
\[
\CW(\pi_f) := a_* \ujast {}^{\ua} \CV' \cdot e_{\pi_f} \; .
\]
\end{Def}

Definition~\ref{5E} should be compared to \cite[Sect.~4.2.0]{Sc}.
Given our construction, the following is obvious.

\begin{Thm} \label{5F}
Assume $\ua = \alpha(k_1,k_2,r)$ to be regular, 
\emph{i.e.}, $k_1 > k_2 > 0$. Let $L$ be a field of characteristic zero,
and $Y_{\pi_f}$ a minimal right ideal
of $R(\FH)\otimes_F L$.
The realizations of the motive $\CW(\pi_f)$ associated to $Y_{\pi_f}$
are concentrated in the single cohomological degree $r+3$, 
and they take the values
$W(\pi_f)$
(in the Hodge theoretic setting) resp.\ 
$W(\pi_f)_\ell$ (in the $\ell$-adic setting).
\end{Thm}

A special case occurs when $Y_{\pi_f}$ is of dimension one over $L$, 
\emph{i.e.}, corres\-ponds to a non-trivial character of $R(\FH)$ with
values in $L$. The automorphic form is then an eigenform for the Hecke
algebra. This is the analogue of 
the situation considered in \cite{Sc} for elliptic
cusp forms. \\

The motive $\CW(\pi_f)$ being a direct factor of $a_* \ujast {}^{\ua} \CV'$,
our results on the latter from Section~\ref{3} have obvious consequences
for the realizations of $\CW(\pi_f)$.

\begin{Cor}  \label{5G}
Assume $\ua = \alpha(k_1,k_2,r)$ to be regular. Let $L$ be a field of characteristic zero,
and $Y_{\pi_f}$ a minimal right ideal
of $R(\FH)\otimes_\BQ L$. Let $p$ be a prime number not dividing 
the level of $K$. 
Let $\ell$ be different from $p$. \\[0.1cm]
(a)~The $p$-adic realization $W(\pi_f)_p$ of 
$\CW(\pi_f)$ is crystalline. \\[0.1cm]
(b)~The $\ell$-adic realization $W(\pi_f)_\ell$ of $\CW(\pi_f)$
is unramified at $p$. \\[0.1cm]
(c)~The characteristic polynomials of the following coincide: (1)~the action
of Frobenius $\phi$ on the $\phi$-filtered module associated to $W(\pi_f)_p$, 
(2)~the action of a geometrical Frobenius automorphism at
$p$ on $W(\pi_f)_\ell$.
\end{Cor}

\begin{Proof}
Parts~(a) and (b) follow from Remark~\ref{3H}.

As for (c), in order to apply \cite[Thm.~2.~2)]{KM},
use that both realizations are cut out by the \emph{same}
cycle from the cohomology of a smooth and proper scheme
over the field $\BF_p$ 
(cmp.\ the proof of Corollary~\ref{3I}).
\end{Proof}

Corollary~\ref{5G} should be compared to \cite[Thm.~1.2.4]{Sc}.

\begin{Rem}
Part~(c) of Corollary~\ref{5G}
is already contained in \cite[Thm.~1]{U}.
\end{Rem}


\bigskip

%
%


\begin{thebibliography}{99}

\bibitem[A]{Anc}
G.~Ancona, {\it D\'ecomposition de motifs ab\'eliens}, 
manuscripta math.~{\bf 146} (2015), 307--328.

\bibitem[BBD]{BBD}
A.A.~Beilinson, J.~Bernstein, P.~Deligne,
{\it Faisceaux pervers},
in B.~Teissier, J.L.~Verdier (eds.),
{\it Analyse et topologie sur les espaces singuliers (I)},
Ast\'erisque~{\bf 100},
Soc.\ Math.\ France (1982).

\bibitem[Be]{Be}
A.A.~Beilinson,
{\it On the crystalline period map},
Camb.~J.~Math.~{\bf 1} (2013), 1--51.

\bibitem[Bo]{Bo}
M.V.~Bondarko,
{\it Weight structures vs.\ $t$-structures; weight filtrations, 
spectral sequences, and complexes (for motives and in general)},
J.~$K$-Theory~{\bf 6} (2010), 387--504.

\bibitem[BW]{BW}
J.I.~Burgos, J.~Wildeshaus,
{\it Hodge modules on Shimura varieties and their higher
direct images in the Baily--Borel compactification},
Ann.\ Scient.\ ENS~{\bf 37} (2004), 363--413.

\bibitem[CD]{CD}
D.-C.~Cisinski, F.~D\'eglise,
{\it Triangulated categories of mixed motives},
preprint, December~2009, version dated December~18, 2012, 279~pages, 
available on arXiv.org under 
{\tt http://arxiv.org/abs/0912.2110} 

\bibitem[C]{C}
G.~Clo\^{\i}tre,
{\it On the interior motive of certain Shimura varieties: the case of Picard varieties}, 
preprint, May~2017, 32~pages, 
available on arXiv.org under 
{\tt https://arxiv.org/abs/1705.03235}

\bibitem[DN]{DN}
F.~D\'eglise, W.~Nizio\l, 
{\it On $p$-adic absolute Hodge cohomology and syntomic coefficients. I},
Comment.\ Math.\ Helv.~{\bf 93} (2018), 71--131. 

\bibitem[DM]{DM}
C.~Deninger, J.P.~Murre, 
{\it Motivic decomposition of abelian schemes and the Fourier transform},
J.\ reine angew.\ Math.~{\bf 422} (1991), 201--219. 

\bibitem[FC]{FC}
G.~Faltings, C.-L.~Chai 
{\it Degeneration of abelian varieties},
Erg.\ der Mathematik und ihrer Grenzgeb.\ 3.~Folge~{\bf 22},
Springer--Verlag Berlin (1990).

\bibitem[Ha]{Ha}
G.~Harder,
{\it Cohomology of Arithmetic Groups},
book in preparation, version dated December 12, 2017, 329~pages,
available under
{\tt http://www.math.uni-bonn.de/people/harder/Manuscripts/buch/
Volume-III.pdf}

\bibitem[H\'e]{H}
D.~H\'ebert,
{\it Structures de poids \`a la Bondarko sur les motifs de Beilinson},
Compositio Math.~{\bf 147} (2011), 1447--1462.

\bibitem[KM]{KM}
N.M.~Katz, W.~Messing,
{\it Some Consequences of the Riemann Hypo\-the\-sis for Varieties 
over Finite Fields},
Invent.~Math.~{\bf 23} (1974), 73--77.

\bibitem[K]{Ki}
S.-I.~Kimura,
{\it Chow groups are finite-dimensional, in some sense},
Math.\ Ann.~{\bf 331} (2005), 173--201.

\bibitem[Lm]{Lem}
F.~Lemma,
{\it On higher regulators of Siegel threefolds I: The vanishing on the boundary},
Asian J.\ Math.~{\bf 19} (2015), 83--120.

\bibitem[Lv1]{L2}
M.~Levine,
{\it Smooth Motives},
in R.~de Jeu, J.D.~Lewis (eds.),
{\it Motives and Algebraic Cycles. A Celebration in Honour of Spencer J.~Bloch},
Fields Institute Communications~{\bf 56}, 
Amer.\ Math.\ Soc.\ (2009), 175--231.

\bibitem[Lv2]{L3}
M.~Levine,
{\it Tate motives and the fundamental group},
in V.~Srinivas (ed.),
{\it Cycles, Motives and Shimura Varieties},
Tata Institute of Fundamental Research (2010),
265--392.

\bibitem[MT]{MT}
A.~Mokrane, J.~Tilouine,
{\it Cohomology of Siegel varieties with $p$-adic integral coefficients and applications},
in A.~Mokrane, P.~Polo, J.~Tilouine (eds.),
{\it Cohomology of Siegel varieties},
Ast\'erisque~{\bf 280},
Soc.\ Math.\ France (2002), 1--95.

\bibitem[M]{M}
M.V.~Moya Giusti,
{\it Ghost classes in the cohomology of the Shimura varie\-ty associated to $GSp_4$}, 
Proc.\ Amer.\ Math.\ Soc.~{\bf 146} (2018), 2315--2325.

\bibitem[O'S]{O'S}
P.~O'Sullivan,
{\it Algebraic cycles on an abelian variety},
J.\ Reine Angew.\ Math.~{\bf 654} (2011), 1--81.

\bibitem[P1]{P}
R.~Pink,
{\it Arithmetical compactification of mixed Shimura varieties},
Bonner Mathematische Schriften~{\bf 209}, Univ.\ Bonn (1990).

\bibitem[P2]{P2}
R.~Pink,
{\it On $\ell$-adic sheaves on Shimura varieties and their higher direct
images in the Baily--Borel compactification},
Math.\ Ann.~{\bf 292} (1992), 197--240.

\bibitem[Sap]{Sa}
L.~Saper, 
{\it $\CL$-modules and the conjecture of Rapoport and Goresky--MacPherson}, 
in J.~Tilouine, H.~Carayol, M.~Harris, M.-F.~Vign\'eras (eds.), 
{\it Automorphic forms. I}, 
Ast\'erisque~{\bf 298} (2005), 319--334.

\bibitem[Sat]{Sat}
I.~Satake, 
{\it Compactification des espaces quotients de Siegel, II}, 
S\'eminaire Henri Cartan~{\bf 10} (1957--58), exp.~no.~13.

\bibitem[Sc]{Sc}
A.J.~Scholl,
{\it Motives for modular forms},
Invent.\ Math.~{\bf 100} (1990), 419--430.

\bibitem[Sh]{Sh}
G.~Shimura,
{\it Introduction to the arithmetic theory of automorphic functions},
Publications of the Math.\ Soc.\ of Japan.~{\bf 11}, Iwanami Shoten, Publishers and Princeton Univ.\ Press (1971).

\bibitem[U]{U}
E.~Urban,
{\it Sur les repr\'esentations $p$-adiques associ\'ees aux repr\'esentations cuspidales de $GSp_{4/\BQ}$},
in J.~Tilouine, H.~Carayol, M.~Harris, M.-F.~Vign\'eras (eds.),
{\it Formes automorphes (II). Le cas du groupe $GSp(4)$},
Ast\'erisque~{\bf 302},
Soc.\ Math.\ France (2005), 151--176.

\bibitem[V]{Vo}
D.A.~Vogan,
{\it Representations of Real Reductive Lie Groups},
Prog.\ in Math.~{\bf 15},
Birkh\"auser (1981).

\bibitem[We]{We}
R.~Weissauer,
{\it Four dimensional Galois representations},
in J.~Tilouine, H.~Carayol, M.~Harris, M.-F.~Vign\'eras (eds.),
{\it Formes automorphes (II). Le cas du groupe $GSp(4)$},
Ast\'erisque~{\bf 302},
Soc.\ Math.\ France (2005), 67--150.

\bibitem[Wi1]{W1}
J.~Wildeshaus,
{\it The canonical construction of mixed sheaves on mixed Shimura varieties},
in: {\it Realizations of Polylogarithms},
Lect.\ Notes Math.~{\bf 1650},
Springer-Verlag (1997), 77--140.

\bibitem[Wi2]{W3}
J.~Wildeshaus,
{\it On the boundary motive of a Shimura variety},
Compositio Math.~{\bf 143} (2007), 959--985.

\bibitem[Wi3]{W4}
J.~Wildeshaus,
{\it Chow motives without projectivity},
Compositio Math.~{\bf 145} (2009), 1196--1226.

\bibitem[Wi4]{W10}
J.~Wildeshaus,
{\it Motivic intersection complex},
in J.I.~Burgos and J.~Lewis (eds.), 
{\it Regulators~III. Proceedings of the conference held
at the University of Barcelona, July 11--23, 2010},
Contemp.\ Math.~{\bf 571}, Amer.\ Math.\ Soc.\ (2012), 255--276.

\bibitem[Wi5]{W7}
J.~Wildeshaus,
{\it On the Interior Motive of Certain Shimura Varieties: 
the Case of Hilbert--Blumenthal varieties},
Int.\ Math.\ Res.\ Notices~{\bf 2012} (2012), 2321--2355.

\bibitem[Wi6]{W8}
J.~Wildeshaus,
{\it On the interior motive of certain Shimura varieties: the case of Picard surfaces}, 
manuscripta math.~{\bf 148} (2015), 351--377.

\bibitem[Wi7]{W9}
J.~Wildeshaus,
{\it Intermediate extension of Chow motives of Abelian type}, 
Adv.\ Math.~{\bf 305} (2017), 515--600.

\bibitem[Wi8]{W11}
J.~Wildeshaus,
{\it Weights and conservativity}, 
Algebraic Geometry~{\bf 5} (2018), 686--702.

\bibitem[Wi9]{W12}
J.~Wildeshaus,
{\it Chow motives without projectivity, II}, 
to appear in Int.\ Math.\ Res.\ Notices, DOI: {\tt 10.1093/imrn/rny257},
available on arXiv.org under 
{\tt http://arxiv.org/abs/1705.10502}

\end{thebibliography}
\end{document}